\documentclass[10pt]{amsart}
\usepackage[margin=2.5cm]{geometry}

\usepackage[T1]{fontenc}
\usepackage[utf8]{inputenc}
\usepackage{lmodern}
\usepackage[british]{babel}

\usepackage{csquotes}
\usepackage[
    sorting=nyt,
    doi=true,
    url=false,
    maxbibnames=4,
    backend=biber,
    giveninits=true,
    style=numeric-comp
    ]{biblatex}
\bibliography{literature}

\usepackage[dvipsnames]{xcolor}
  \definecolor{darkred}{RGB}{139,0,0}
  \definecolor{mediumblue}{RGB}{0,0,205}
  \definecolor{forestgreen}{RGB}{34,139,34}
\usepackage{graphicx}

\usepackage[colorlinks=true]{hyperref}
\hypersetup{linkcolor=darkred, urlcolor=forestgreen, citecolor=mediumblue}
\usepackage{cleveref}

\usepackage{amsmath, amsfonts, amssymb}
\usepackage{mathtools}
\usepackage{commath}

\usepackage{bm}
\usepackage{calrsfs}
\usepackage{stmaryrd}
\DeclareMathAlphabet{\pazocal}{OMS}{zplm}{m}{n}
\DeclareMathAlphabet{\pazocalbf}{OMS}{cmsy}{b}{n}
\usepackage{etoolbox}
\usepackage[normalem]{ulem}
\preto\subequations{\ifhmode\unskip\fi}
\newcommand{\clem}[2][]{\ifstrempty{#1}{%
  c_{\hyperref[#2]{L\ref*{#2}}}}{%
  c_{\hyperref[#2]{L\ref*{#2}#1}}}}

\newcommand{\ceq}[1]{c_{\hyperref[#1]{(\ref*{#1})}}}

\usepackage{amsthm}
\theoremstyle{plain}

\newtheorem{lemma}{Lemma}

\theoremstyle{definition}

\newtheorem{assumption}{Assumption}

\theoremstyle{remark}
\newtheorem{remark}{Remark}

\usepackage{booktabs}
  \AtBeginEnvironment{tabular}{\footnotesize}
\usepackage{siunitx}
\newcommand{\tend}{T}
\renewcommand{\O}{\Omega}
\newcommand{\G}{\Gamma}

\let\div\relax
\DeclareMathOperator{\div}{div}
\DeclareMathOperator{\dist}{dist}
\DeclareMathOperator{\meas}{meas}
\DeclareMathOperator{\Int}{Int}
\DeclareMathOperator{\sign}{sign}

\newcommand{\wb}{\bm{w}}
\newcommand{\bw}{\wb}
\newcommand{\nb}{\bm{n}}
\newcommand{\xb}{\bm{x}}
\newcommand{\bn}{\bm{n}}
\newcommand{\ahn}{a_h^n}

\newcommand{\Sc}{\pazocal{S}}
\newcommand{\Sd}{\Sc_{\delta}}
\newcommand{\Sdh}{\Sc_{\delta_h}^\pm}
\newcommand{\Sdhp}{\Sc_{\delta_h}^+}

\newcommand{\Ec}{\pazocal{E}}

\newcommand{\dt}{\Delta t}
\newcommand{\OO}{\pazocal{O}}

\newcommand{\Odh}{\OO_{\delta_h}}
\newcommand{\Tc}{\pazocal{T}}
\newcommand{\Th}{\Tc_h}
\newcommand{\Thn}{\Th^n}
\newcommand{\ThS}{\Tc_{\Sc^{\pm}_{\delta}}^n}
\newcommand{\ThSp}{\Tc_{\Sc^{+}_{\delta}}^n}

\newcommand{\Fc}{\pazocal{F}}
\newcommand{\Fh}{\Fc_h}
\newcommand{\Fhn}{\Fh^n}
\newcommand{\OdhT}{\widetilde{\O}^n_h}
\newcommand{\OdhTlast}{\widetilde{\O}^{n-1}_h}
\newcommand{\OdhTone}{\widetilde{\O}^1_h}

\newcommand{\Oh}{\O_{h}}
\newcommand{\On}{\O^{n}}
\newcommand{\Ohn}{\Oh^{n}}
\newcommand{\uh}{u_{h}}
\newcommand{\uhn}{\uh^{n}}
\newcommand{\vh}{v_h}
\newcommand{\wh}{w_h}
\newcommand{\Gn}{\G^{n}}
\newcommand{\Ghn}{\Gn_{h}}

\newcommand{\RR}{\mathbb{R}}
\newcommand{\PP}{\mathbb{P}}
\newcommand{\NN}{\mathbb{N}}

\newcommand{\jump}[1]{\llbracket#1\rrbracket}
\newcommand{\restr}[2]{{\left.\kern-\nulldelimiterspace#1\right|_{#2}}}

\renewcommand{\norm}[2]{\Vert #1 \Vert_{#2}}
\newcommand{\eps}{\varepsilon}

\newcommand{\eoc}[1]{\text{eoc}_{#1}}

\title[A conservative Eulerian FEM for transport and diffusion in moving domains]{A conservative Eulerian finite element method for transport and diffusion in moving domains}

\author[M. Olshanskii]{Maxim Olshanskii$^{\ast}$}\thanks{$^{\ast}$ ORCID: 0000-0002-9102-6833}
\address{Department of Mathematics, University of Houston, 651 PGH Houston, Texas 77204, USA}
\email{molshan@math.uh.edu}

\author[H. v. Wahl]{Henry von Wahl$^{\dagger}$}\thanks{$^{\dagger}$ ORCID: 0000-0002-0793-1647}
\address{Friedrich-Schiller-Universität, Fakultät für Mathematik und Informatik, 
  Ernst-Abber-Platz 2, 07743 Jena, Germany}
\email{henry.von.wahl@uni-jena.de}
\address{Institute for Computational and Experimental Research in Mathematics, 
  Brown University, 121 South Main Street, Box E, Providence, RI 02903, USA}

\date{\today}

\subjclass{65M12, 65M60, 65M85}

\keywords{Eulerian time stepping, evolving domains, discrete conservation, ghost penalty}

\begin{document}

\begin{abstract}
The paper introduces a finite element method for an Eulerian formulation of partial differential equations governing the transport and diffusion of a scalar quantity in a time-dependent domain. The method follows the idea from Lehrenfeld \& Olshanskii [ESAIM: M2AN, 53(2): 585-614, 2019] of a solution extension to realise the Eulerian time-stepping scheme. However, a reformulation of the partial differential equation is suggested to derive a scheme which conserves the quantity under consideration exactly on the discrete level. For the spatial discretisation, the paper considers an unfitted finite element method. Ghost-penalty stabilisation is used to realise the discrete solution extension and gives a scheme robust against arbitrary intersections between the mesh and geometry interface. The stability is analysed for both first- and second-order backward differentiation formula versions of the scheme. Several numerical examples in two and three spatial dimensions are included to illustrate the potential of this method.
\end{abstract}

\maketitle

\section{Introduction}
\label{sec:intro}

Mathematical models from biology, chemistry, physics and engineering often include partial differential equations (PDEs) posed on evolving domains. Examples of such problems are fluid flows in or around moving structures such as wind turbines~\cite{PBS19} or particle-laden flows~\cite{Sun03}, biomedical applications such as blood flow through the cardiovascular system~\cite{vassilevski2020personalized} or tumour growth~\cite{GLNS18} and multi-phase flows such as rising bubbles~\cite{JB17} or droplets in microfluidic devices~\cite{CK19}.

Computational methods for such problems face the challenge of dealing with both inherently Eulerian quantities, such as concentrations and temperature, and inherently Lagrangian quantities, such as displacements. Many methods that deal with PDEs on evolving domains are based on the Lagrangian or arbitrary Lagrangian–Eulerian (ALE) formulations~\cite{hirt1974arbitrary}. Lagrangian and ALE methods can be  based on a reference configuration, into which the problem is mapped and in which it is then solved using either standard time-stepping schemes or space-time Galerkin formulations; see, e.g.,~\cite{tezduyar1992new,masud1997space}. An advantage of this approach is that fitted and adapted meshes of the reference geometry can be used, leading to good resolution of the moving interface. However, the approach becomes more limited when the deformation gets large and fails in the case of topology changes. To avoid these issues, purely Eulerian methods can be considered. This will also be our approach.

Within the framework of finite element methods (FEM), unfitted finite element methods have gained traction in the past decade. These include the eXtenden finite element method (XFEM)~\cite{fries2010extended}, CutFEM~\cite{BCH+14}, the finite cell method~\cite{parvizian2007finite} and TraceFEM~\cite{olshanskii2009finite}. These methods are of particular interest in situations of complex geometries, where mesh generation can be a difficult task, as they separate the mesh from the geometry description. As such, they are well suited to moving domain problems, where complex geometries can easily arise. While these methods have been extensively studied for problems on stationary domains, unfitted finite element methods for moving domain problems are less established. One difficulty in the context of evolving domains is due to the approximation of the Eulerian time-derivative
\begin{equation*}
  \partial_t u \approx \frac{1}{\dt}(u(t^n) - u(t^{n-1})).
\end{equation*}
When the domain of interest $\Omega$ changes in time, then $u(t^n)$ is defined on $\O(t^n)\subset\RR^d$ and $u(t^{n-1})$ on $\O(t^{n-1})$. Consequently, the difference $u(t^n) - u(t^{n-1})$ may no longer be well-defined.

Several approaches have been developed to try and remedy this problem. These include characteristic-based approaches~\cite{HP17}, applying the ALE approach within one time step and projecting the solution back into the original reference configuration~\cite{CHCB09} and a Galerkin time-discretisation using modified quadrature rules to recover a classical time-stepping scheme~\cite{FR17}. However, these approaches require expensive projections between domains at different timers. A different and widely studied approach is to use space-time Galerkin approaches~\cite{LR13,ORX14,HLZ16,Pre18,AB21,HLP23}. Space-time methods require the solution of $d+1$-dimensional problems and are consequently computationally expensive.

In this work, we will focus on CutFEM as our spatial discretisation. To this end, we consider a static background mesh of our computational domain, and we assume that the volume domain of interest $\O(t)$ evolves smoothly within this static computational domain. To handle instabilities arising due to arbitrary `cuts' between $\O(t)$ and the background mesh, CutFEM uses ghost-penalty stabilisation~\cite{Bur10}. For the temporal discretisation, we modify the Eulerian time-stepping scheme suggested in~\cite{LO19}, based on implicit extensions of the solution to a neighbourhood of order $\OO(\dt)$ using ghost-penalty stabilisation.

The idea of using the extension provided by ghost-penalty stabilisation in CutFEM to enable an Eulerian time-stepping scheme for moving domain problems can be traced back to~\cite{Sch17}. However, the method there was limited to a geometric CFL condition $\dt \leq c\,h$, as the extension was computed in a separate step. Furthermore, the method was not analysed. The more general method in~\cite{LO19} was then extended and analysed in a number of different settings. In~\cite{LL21}, the method was expanded to higher-order geometry approximations using isoparametric mappings. The method was extended to the Stokes problem on moving domains in~\cite{BFM22} and~\cite{vWRL20} using stabilised equal order and Taylor-Hood elements, respectively. In~\cite{NO24}, the approach was extended to the linearised Navier-Stokes problem in evolving domains and the a priori error analysis was improved upon. Furthermore, in~\cite{vWR21a,vWR23}, the approach was analysed for a parabolic model problem for coupled fluid-rigid body interactions. All these approaches were based on a backward differentiation formula (BDF) discretisation of the time-derivative. In~\cite{FS22}, the approach has been extended to a Crank-Nicolson discretisation of the time-derivative for the heat-equation on moving domains. However, while the latter study is the first work to extend the approach to non-BDF time-stepping schemes, the analysis required a parabolic CFL condition $\dt\leq c\,h^2$ not needed for the BDF-type scheme.

In this paper, we revisit the equations governing the conservation of a scalar quantity with a diffusive flux in a domain passively evolved by a smooth flow, as in the original work~\cite{LO19}. While optimal order of convergence with respect to the spatial, temporal- and geometry-approximation was proved in that paper, the method was not optimal in the sense that the conservation of total mass can be violated by the discrete solution. The focus of this work will, therefore, be to develop a modification of this Eulerian time-stepping scheme to preserve the mass of the scalar variable exactly on the discrete level.

The main idea of our approach is to rewrite the PDE problem using an identity derived from the Reynolds transport theorem. Using this identity as the basis of our finite element method, we arrive at a scheme for which we can show that the discrete solution conserves the total mass. The stability analysis presented is similar to that of~\cite{LO19}. 

The remainder of this paper is structured as follows. In \Cref{sec:problem}, we briefly present the PDE problem under consideration, cover the temporally semi-discrete version of our scheme and analyse the unique solvability of the problem and the stability of the scheme. In this semi-discrete setting, the method applies the extension to the test function in the problem formulation, whereas the method in~\cite{LO19} and subsequent works use an extension of the solution into the next domain. In \Cref{sec:discrete}, we then cover the fully discrete method. We present the fully discrete scheme based on both BDF1 and BDF2 formulas to discretise the time-derivative in the scheme. Furthermore, we consider the discrete stability of both schemes. Finally, in \Cref{sec:num-ex}, we consider a number of numerical examples in two and three spatial dimensions. In examples with a given analytical solution, we investigate the numerical convergence of the schemes with respect to mesh refinement, time-step refinement and combined mesh-time step refinement.

We finally note that a rigorous error analysis of the method remains an open problem. Numerical experiments demonstrate optimal convergence rates in some examples and slightly sub-optimal rates in certain cases as demonstrated in \Cref{sec:num-ex}. The latter observation warrants further investigation of conservative unfitted FEMs.

\section{Problem setup and preliminaries}
\label{sec:problem}
Let us consider a time-dependent domain $\O(t)\subset\RR^d$, with $d\in\{2,3\}$, which  evolves smoothly and for all time  $t\in[0, \tend]$, $\tend>0$, stays bounded and embedded in an background fixed domain $\widetilde\O$,  $\O(t)\subset\widetilde\O\subset\RR^d$.
We further assume that the evolution of $\O(t)$ is a passive transport by a smooth flow field $\bw:\widetilde\O\times(0, \tend)\to \RR^d$, i.e., the Lagrangian mapping $\Psi(t)\,:\,\O_0\to\O(t)$ is given by  $\Psi(t,y)$ that solves the ODE system
\begin{equation}\label{Lagrange}
\Psi(0,y)=y,\quad  \frac{\partial \Psi(t,y)}{\partial t}=\bw(t,\Psi(t,y)),\quad t\in[0,T],\quad y\in\O(0),
\end{equation}
where we assume the reference domain $\O_0$ is piecewise smooth and Lipschitz.

The conservation of a scalar quantity with a diffusive flux in  $\O(t)$ and no mass exchange through the boundary is governed by the equations
\begin{subequations}\label{eqn.strong-problem}
\begin{align}
  && \partial_t u + \div(u\wb) - \nu \Delta u &= f &&\text{in }\O(t),\label{eqn.strong-problem.1}\\
  && \nabla u \cdot \nb &= 0 &&\text{on }\G(t),\label{eqn.strong-problem.2}
\end{align}
\end{subequations}
where  $\nu>0$ is a constant diffusion coefficient, $\G(t)\coloneqq \partial\O(t)$ is the boundary of the domain,  and $\nb$ is the unit normal vector on $\G(t)$. The well-posedness of \eqref{eqn.strong-problem} is addressed, for example, in~\cite{alphonse2015abstract}.

The quantity $u$ is globally conserved up to the total contribution of the source term $f$. This is easy to see using Reynolds' transport formula, equation \eqref{eqn.strong-problem.1} and the no-flux boundary condition  \eqref{eqn.strong-problem.2}:
\begin{equation}\label{eq:balance}
\begin{split}
\frac{\dif}{\dif t} \int_{\Omega(t)} u \dif{x} &=  \int_{\Omega(t)} \frac{\partial}{\partial t} u \dif{x} + \int_{\partial \Omega(t)} (\bw \cdot \bn) u \dif{s}
=  \int_{\Omega(t)} (\frac{\partial u}{\partial t} + \div(u \bw)\, )\dif{x}
=  \int_{\Omega(t)} (\nu \Delta u +f)\, \dif{x}\\
&=  \int_{\partial \Omega(t)} \nu \nabla u \cdot \bn \dif{s}+\int_{\Omega(t)} f\, \dif{x}=\int_{\Omega(t)} f\, \dif{x}.
\end{split}
\end{equation}
This is a fundamental property, which we would like to preserve in a numerical method. 

More generally, for  any $v\in H^1(\widetilde\O)$, we have 
\begin{equation*}
  \frac{\dif}{\dif t} \int_{\O(t)} u v \dif x = \int_{\O(t)} \partial_t uv \dif x + \int_{\G(t)} \wb\cdot\nb uv \dif s
  = \int_{\O(t)} ((\partial_t u + \div(u\wb))v + u\wb\cdot\nabla v) \dif x,
\end{equation*}
which follows from Reynolds' transport formula and the observation that $v$ is time-independent.  
This yields the identity
\begin{equation}\label{eqn.identity1}
   \int_{\O(t)} (\partial_t u + \div(u\wb))v \dif x = \frac{\dif}{\dif t}\int_{\O(t)} uv \dif x - \int_{\O(t)} u\wb\cdot\nabla v \dif x.
\end{equation}
Therefore, multiplying \eqref{eqn.strong-problem.1} with a test-function $v\in H^1(\widetilde\O)$, integrating by parts and using \eqref{eqn.strong-problem.2} and \eqref{eqn.identity1} gives the following identity satisfied for a smooth solution to \eqref{eqn.strong-problem} and any $v\in H^1(\widetilde\O)$
\begin{equation}\label{eqn.weak-problem}
  \frac{\dif}{\dif t}\int_{\O(t)} uv \dif x + \int_{\O(t)} (\nu\nabla u(t) - u(t)\wb(t))\cdot\nabla v\dif x = \int_{\O(t)}f(t)v\dif x\quad\text{for}~t\in[0,\tend].
\end{equation}
 This identity will form the basis of our discretisation.

\subsection{Temporal semi-discretisation} To  elucidate our numerical method construction,  we first formulate a semi-discrete method. For this, we need a smooth extension operator
\begin{equation*}
  \Ec: L^2(\O(t)) \rightarrow L^2(\widetilde\O),
\end{equation*}
such that the function remains unchanged in the original domain $\O(t)$ and it holds that
\begin{equation*}
  \norm{\Ec u}{H^k(\widetilde\O)} \lesssim \norm{u}{H^k(\O(t))}
  \quad\text{and}\quad
  \norm{\nabla \Ec u}{\widetilde\O} \lesssim \norm{\nabla u}{\O(t)}.
\end{equation*}
We refer to~\cite{LO19} for the details of the explicit construction of such an extension operator based on the classical linear and continuous universal extension operator for Sobolev spaces from~\cite{Ste70}.

Let us consider a uniform time step $\dt = T / N$ for a fixed $N\in\NN$. We denote $t_n = n\dt$, $\On = \O(t_n)$ and $u^n \approx u(t^n)$, and define
\begin{equation*}
  a^n(u,v) \coloneqq \int_{\On}\nu\nabla u\cdot\nabla v - u (\wb^n\cdot\nabla v) \dif x
  \quad\text{and}\quad
  f^n(v) \coloneqq \int_{\On}f(t^n) v \dif x.
\end{equation*}
Using an implicit Euler (or BDF1) discretisation in time of \eqref{eqn.weak-problem} then leads to the problem: Given $u^{n-1}\in H^1(\O^{n-1})$, 
find $u^{n}\in H^1(\On)$ solving 
\begin{equation}\label{eqn:scheme.bdf1.smooth}
  \frac{1}{\dt}\left[\int_{\On} u^nv \dif x - \int_{\O^{n-1}} u^{n - 1}\Ec v \dif x\right] + a^n(u^n,v) = f^n(v),
\end{equation}
for all $v\in H^1(\On)$. Note that unlike~\cite{LO19} and subsequent studies, the extension operator in the semi-discrete method acts on the test function rather than on the solution.  
\medskip

In what follow, we shall denote the inner-product for functions in $L^2(S)$, for some domain $S$ as $(\cdot, \cdot)_S$. Similarly, $\norm{\cdot}{S}$ will denote the $L^2(S)$ norm. Furthermore, we use the notation $a\lesssim b$, iff there exists a constant $c>0$, independent of $\dt$, $h$ and the mesh-interface cut position, such that $a\leq cb$. Similarly, we write $a\gtrsim b$, iff $a\geq cb$, and $a\simeq b$, iff both $a\lesssim b$ and $a\gtrsim b$ hold.

For the bilinear form $a^n(\cdot,\cdot)$, we have the following Gårding's inequality:
\begin{lemma}
Let $w_\infty\coloneqq \max_{t\in[0,T]}\norm{\wb(t)}{\infty}$. It holds that
\begin{equation*}
  a^n(u,u) \geq \frac{\nu}{2}\Vert\nabla u\Vert_{\On}^2 - \frac{w_\infty^2}{\nu}\Vert u\Vert_{\On}^2,
\end{equation*}
i.e., we have unique solvability of the problem \eqref{eqn:scheme.bdf1.smooth} in every time step, if
\begin{equation*}
  \dt < \nu w_\infty^{-2}. 
\end{equation*} 
\end{lemma}

\begin{proof}
Follows by applying the Cauchy-Schwarz inequality and a weighted Young's inequality.
\end{proof}

\begin{remark}
The above time step can be limiting. In practice, we have not found such a time-step restriction to be necessary. Indeed, one can show that the problem is uniquely solvable for $\dt> \nu w_\infty^{-2}$, if $1 / \dt$ does not belong to the spectrum of the operator $A$ corresponding to the bilinear form $a^n(\cdot,\cdot)$ \cite{BS08,Agm65}.
\end{remark}

It would also be instructive to check the stability of the semi-discrete scheme.
We have the following stability result.

\begin{lemma}\label{lemma.smooth-euler.stability}
  The solution of the semi-discrete scheme in \eqref{eqn:scheme.bdf1.smooth}, fulfils the stability estimate
  \begin{equation*}
    \norm{u^N}{\O^N}^2 + \dt\sum_{n=1}^N\frac{\nu}{2}\norm{\nabla u^n}{\On}^2
    \leq \exp(\clem{lemma.smooth-euler.stability}T)\left[\norm{u^0}{\O^0}^2 + \dt\sum_{n=1}^{N}\frac{2}{\nu}\norm{f^n}{H^{-1}(\On)}^2 \right],
  \end{equation*}
  with a constant $\clem{lemma.smooth-euler.stability}>0$ independent of $\dt$, $N$ and $u$, once $\dt$ is sufficiently small.
\end{lemma}

\begin{proof}
Define the $\delta$-strip around the moving domain boundary as
\begin{equation*}
  \Sd^{\pm}(\O(t)) \coloneqq \{\xb\in\RR^d : \dist(\xb, \G(t)) \leq \delta\},
\end{equation*}
and the  $\delta$-strip outside of the physical domain as
\begin{equation*}
  \Sd^+(\O(t)) \coloneqq \Sd^{\pm}(\O(t)) \setminus \O(t).
\end{equation*}
We choose $\delta\simeq\Delta t$ but sufficiently large such that $\O^{n-1}\subset \On\cup\Sd^{\pm}(\On)$ for all $n$.
Testing \eqref{eqn:scheme.bdf1.smooth} with an appropriate extension of the solution $v=2\dt\Ec u^n$ leads to
\begin{multline*}
  2\norm{u^n}{\On}^2 - \norm{\Ec u^{n}}{\O^{n-1}}^2 - \norm{u^{n-1}}{\O^{n-1}}^2 + \norm{\Ec u^n - u^{n-1}}{\O^{n-1}}^2 + 2\dt\nu\norm{\nabla u^n}{\On}^n - 2\dt(u^n,\wb\cdot\nabla u^n)_{\On}\\
   = 2\dt(f^n,u^n)_{\On}.
\end{multline*}
Therefore, using the Cauchy-Schwarz inequalities, we have
\begin{multline*}
  \norm{u^n}{\On}^2 + \norm{u^n}{\On\setminus\O^{n-1}}^2
  + \norm{\Ec u^n - u^{n-1}}{\O^{n-1}}^2 + 2\dt\nu\norm{\nabla u^n}{\On}^2\\
  \leq 
    \norm{u^{n-1}}{\O^{n-1}}^2 
    + 2\dt\norm{f^n}{H^{-1}(\On)}\norm{\nabla u^n}{\On}
    + 2\dt\sqrt{2}w_\infty\norm{u^n}{\On}\norm{\nabla u^n}{\O^n}
    + \norm{\Ec u^n}{\O^{n-1}\setminus\On}^2.
\end{multline*}
With two weighted Young's inequalities, we then get
\begin{equation*}
  \norm{u^n}{\On}^2 
  + \dt\nu\norm{\nabla u^n}{\On}^2\\
  \leq 
    \norm{u^{n-1}}{\O^{n-1}}^2 
    + \frac{2\dt}{\nu}\norm{f^n}{H^{-1}(\On)}^2
    + \frac{4\dt}{\nu}w_\infty^2\norm{u^n}{\On}^2
    + \norm{\Ec u^n}{\O^{n-1}\setminus\On}^2.
\end{equation*}
Since $\O^{n-1}\setminus\On\subset\Sd(\On)$, we have the estimate
\begin{equation}\label{eqn:smooth-extensio-strip-bound}
  \norm{\Ec u^n}{\O^{n-1}\setminus\On}^2 \leq
  \norm{\Ec u^n}{\Sd(\On)}^2
  \lesssim (1 + \eps^{-1})\delta\norm{u^n}{\On}^2 + \delta \eps \norm{\nabla u^n}{\On}^2,
\end{equation}
see \cite[Lemma~3.5]{LO19} for the proof of the second estimate. With the choice of $\delta\simeq\Delta t$, we may set $\eps$, such that
\begin{equation*}
  \norm{u^n}{\On}^2 
  + \frac{\dt\nu}{2}\norm{\nabla u^n}{\On}^2\\
  \leq 
    \norm{u^{n-1}}{\O^{n-1}}^2 
    + \frac{2\dt}{\nu}\norm{f^n}{H^{-1}(\On)}^2
    + \dt\left(\frac{4w_\infty^2}{\nu} + c(\nu^{-1})\right)\norm{u^n}{\On}^2.
\end{equation*}
Summing over $n=1,\dots, N$ and applying a discrete Grönwall's estimate, cf.~\cite{HR90}, for $\dt$ sufficiently small such that $\dt\big(4w_\infty^2\nu^{-1} + c(\nu^{-1})\big) < 1$, gives the result.
\end{proof}

\section{Discrete Method}
\label{sec:discrete}
\begin{figure}
  \centering
  \includegraphics{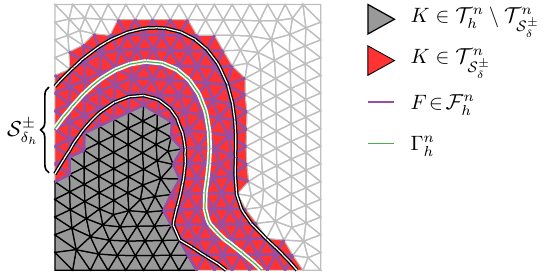}
  \caption{Sketch of the different sets of elements and facets used in the discrete method.}
  \label{fig:elements}
\end{figure}

We now formulate a fully discrete counterpart to \eqref{eqn:scheme.bdf1.smooth}. For the spatial discretisation, we use unfitted finite element, also known as CutFEM~\cite{BCH+14}. We consider a quasi-uniform simplicial mesh $\Th$ of $\widetilde{\O}$ with characteristic element diameter $h>0$. Thus, the mesh is independent of $\O(t)$. On this background mesh, we define the finite element space
\begin{equation*}
  V_h \coloneqq \{v\in C(\widetilde\Omega)  : \restr{v}{K}\in\PP^k(K),\text{ for all }K\in\Th \},\quad k\geq1,
\end{equation*}
where $\PP^k(K)$ is space of polynomials of order $k$ on $K$. In each time step, the domain $\On$ is approximated by $\Ohn$, which is computed by an approximated level set function, see \Cref{sec:discrete.subsec:geo} below. We denote the boundary as $\Ghn\coloneqq\partial\Ohn$. For a the discrete formulation, we now consider at each time step $\delta_h$-extensions of the discrete domains:
\begin{equation}\label{eqn:choiche-discrete-delta}
  \Odh(\Ohn):=\{\xb\in\widetilde\O\,:\,\mbox{dist}(\xb,\Ohn)<\delta_h\},\quad \delta_h = c_{\delta} w_{\nb}^\infty \dt,
\end{equation}
with some $c_\delta\simeq 1$ defined later. 
Using this extension, we define the \emph{active mesh} in each time step as the set of elements with some part in the extended domain
\begin{equation*}
  \Thn \coloneqq \{K\in\Th : \dist(\xb,\Ohn) \leq \delta_h,\; \xb\in K \}.
\end{equation*}
Consequently, we define the \emph{active domain} as
\begin{equation*}
  \OdhT = \{\xb\in \overline{K} : K\in\Thn\}.
\end{equation*}
A sketch of this can be seen in \Cref{fig:elements}. With above choice \eqref{eqn:choiche-discrete-delta}, we then have the essential property for our method, that
\begin{equation*}
  \O^{n-1}_h\subset \OdhT.
\end{equation*} 
Finally, in each time step, we then take the restriction of the finite element space to the set of active elements
\begin{equation*}
  V_h^n \coloneqq \{v\in C(\widetilde\Omega)  : \restr{v}{K}\in\PP^k(T),\text{ for all }K\in\Thn \}.
\end{equation*}

\subsection{Finite element formulation}
\label{sec:discrete.subsec:formulation}

To formulate the discrete problem, let us define the (bi-) linear forms
\begin{equation*}
  \ahn(u,v) = \int_{\Ohn}\nu\nabla u\cdot\nabla v - u (\wb^n\cdot\nabla v) \dif x,
  \quad\text{and}\quad
  f_h^n(v)  = \int_{\Ohn}f^n v\dif x.
\end{equation*}
Furthermore, let $s_h^{n}(u,v)$ be a stabilisation bilinear form which acts in a strip $\Sdh$ around $\Ohn$, such that $\Oh^{n+1}\subset \Ohn \cup \OdhT$. In practice, this is a ghost-penalty operator, see \Cref{sec:discrete.subsec:formulation.subsubsec:ghost-penalty} below.

Using an implicit Euler discretisation of the time-derivative, the fully discrete scheme then reads as follows: Given $u^0$, find $\uhn\in V_h^n$, for $n=1,\dots, N$, such that
\begin{equation}\label{eqn:discrete-scheme.bdf1}
  \frac{1}{\dt}\left[\int_{\Ohn} \uhn \vh \dif x - \int_{\Oh^{n-1}}\uh^{n-1} \vh \dif x \right] + \ahn(\uhn, \vh) + \nu s_h^n(\uhn,\vh) = f_h^n(\vh),
\end{equation}
for all $\vh\in V_h^n$.

\begin{remark}[Conservation]
  The scheme presented in \eqref{eqn:discrete-scheme.bdf1} is \emph{conservative}, if the stabilisation term vanishes for constant test functions. Testing with $\vh\equiv 1$, for which $a_h^n$ and the ghost-penalty term vanishes and summing gives
  \begin{equation*}
    \int_{\Ohn} \uhn \dif x = \int_{\Oh^0}\uh^0 \dif x + \dt\sum_{k=1}^{n}\int_{\Oh^k}f^k\dif x,
  \end{equation*}
  which is the discrete counterpart of \eqref{eq:balance}.
\end{remark}

\subsubsection{Higher order in time}
A simple extension of the above method is to use a second-order backward-differentiation formula (BDF2) for the time derivative. This then reads as follows: Given $u^0$, compute $\uh^1$ using \eqref{eqn:discrete-scheme.bdf1} and find $\uhn\in V_h^n$, for $n=2,\dots, N$, such that
\begin{equation}\label{eqn.discrete-scheme:bdf2}
  \frac{1}{\dt}\left[\frac{3}{2}\int_{\Ohn} \uhn \vh \dif x - 2\int_{\Oh^{n-1}}\uh^{n-1} \vh \dif x + \frac{1}{2}\int_{\Oh^{n-2}}\uh^{n-2} \vh \dif x \right] + \ahn(\uhn, \vh) + \nu s_h^n(\uhn,\vh) = f_h^n(\vh),
\end{equation}
for all $\vh\in V_h$. By increasing the width of the extension strip to by a factor of two in \eqref{eqn:choiche-discrete-delta}, this is then well posed.

\begin{remark}
Same as \eqref{eqn:discrete-scheme.bdf1}, the BDF2 scheme conserves the scalar quantity, up to forcing terms. Analogously, the BDF$k$ finite difference methods for $k > 2$ can be applied to discretise the time derivative while ensuring conservation. The authors of \cite{LL21} noted that the stability analysis of BDF3 and BDF4 for a non-conservative formulation can be conducted similarly to that of BDF2. This is likely also true for the conservative formulation; however, we will not pursue this more technical analysis here.
\end{remark}

\subsubsection{Stabilising bilinear form}
\label{sec:discrete.subsec:formulation.subsubsec:ghost-penalty}

In problem \eqref{eqn:discrete-scheme.bdf1}, the extension into the active domain is realised through the stabilising form $s_h^n(\cdot,\cdot)$. We shall take $s_h^n$ to be a ghost-penalty stabilisation operator~\cite{Bur10}. The original aim for this kind of stabilisation is to prevent stability and conditioning problems occurring in unfitted finite elements in the presence of so-called bad cuts, i.e., elements where $|K\cap\O|$ is very small, by giving a stable solution on the whole of each cut element. Ghost-penalty stabilisation additionally provides us with an implicitly extended solution outside of the physical domain. By applying this stabilisation in a strip around $\Ghn$, as first proposed in~\cite{LO19}, we obtain the necessary discrete extension for our Eulerian time-stepping scheme. To this end, we define the set of boundary strip elements as
\begin{equation*}
  \ThS \coloneqq \{K\in\Thn : \dist(\xb,\Ghn) \leq \delta_h,\; \xb\in K \},
\end{equation*}
and the set of facets of these elements, that are interior to the active mesh
\begin{equation*}
  \Fhn \coloneqq \{F=\overline{K}_1\cap\overline{K}_2 : K_1\in \ThS, K_2\in\Thn, K_1\neq K_2, \meas_{d-1}(F) > 0\}.
\end{equation*}

A number of different versions of the ghost-penalty operator exist, the most common version being based on normal derivative jumps across relevant facets, see for example~\cite{BFM22,BH12,FNZ23,FS22,MLLR14}. We prefer the so-called \emph{direct-version}, introduced in~\cite{Pre18}, since this extends more easily to higher-order spaces. For a more detailed overview, we refer to \cite[Section~4.3]{LO19}. 

To define the direct ghost-penalty operator, let $\Ec^\PP\colon \PP^k(K) \rightarrow \PP^k(\RR^d)$ be the canonical global extension operator of polynomials to $\RR^d$, i.e., for $v\in\PP^k(K)$ it holds that $\restr{\Ec^\PP(v)}{K} = v$. For a facet $F=\overline{K}_1\cap\overline{K}_2$, $K_1\neq K_2$, we define the \emph{facet patch} as $\omega_F = \Int(\overline{K}_1\cup\overline{K}_2)$. The volumetric patch jump of a polynomial $\vh\in\PP^{k}(\Th)$ is then defined via
\begin{equation*}
  \restr{\jump{\vh}}{K_1} = \restr{\vh}{K_i} - \Ec^\PP(\restr{\vh}{K_j}),\quad\text{for }i,j\in\{1,2\}\text{ and }i\neq j.
\end{equation*}
The direct version of the ghost-penalty operator is then defined as
\begin{equation*}
  s_h^n(u_h, \vh) \coloneqq \gamma_s \sum_{F\in\Fhn}\frac{1}{h^2}\int_{\omega_F}\jump{u_h}\jump{\vh}\dif x,
\end{equation*}
with a stabilisation parameter $\gamma_s > 0$. For this operator to provide the necessary stability, the stability parameter needs to scale with number of facets which need to be crossed to reach an element interior to the physical domain from an element in the extension strip~\cite{LO19}. This number $L$ depends on the anisotropy between the mesh size and time step by $L \lesssim (1 + \delta_h / h) \simeq (1 + c\dt / h)$. Accordingly, we set
\begin{equation}\label{eqn.gp-parameter-choice}
  \gamma_s = c_\gamma L,\quad\text{with }c_\gamma>0\text{ independent of $h$ and $\dt$}.
\end{equation}
See also~\cite{LO19} for further details to this choice.

\begin{assumption}\label{assumption:ghost-penalty}
For every cut or extension element
\begin{equation*}
  K\in\ThSp \coloneqq \ThS \setminus \{K\in\Thn\;:\; K\subset \Ohn \},
\end{equation*}
there exists an uncut element $K'\in\ThS\setminus\ThSp$, which can be reached by repeatedly crossing facets in $\Fhn$. We assume that the number of facets that must be crossed to reach $K$ from $K'$ is bounded by a constant $L\lesssim 1 + \frac{\delta_h}{h}$. Furthermore, we assume that each uncut element $K'\in\ThS\setminus\ThSp$ is the end of at most $M$ such paths, with $M$ bounded independent of $h$ and $\dt$.
\end{assumption}

\begin{remark}\label{remak:assume-gp}
 We briefly comment on the validity of \Cref{assumption:ghost-penalty} if the domain boundary is smooth and sufficiently well resolved by the mesh, cf.\ \cite[Remark~5.4]{LO19}. Let $x_K$ be the circumcenter of a cut or extension element $K\in \ThSp$. Furthermore, let $p\colon\Odh(\Ohn)\to\Ghn$ be the closet point projection onto the boundary. We know map $x_K$ towards the interior by $y= x_K + \delta_h (p(x_K) - x_K)$. There is an uncut element $K'\in\ThS\setminus\ThSp$ that contains the point $y$ or that can be reached from $y$ by crossing a finite number of facets, which for fine enough mesh depends only on the minimal angle condition. This construction defines a mapping $B:\ThSp\to \ThS\setminus\ThSp$. Now, due to the assumed shape regularity of the mesh, the number of facets $L$ crossed by the path $\{x_K + s (p(x_K) - x_K) \,:\, s\in[0, \delta_h] \}$ is bounded by $c(h + \delta_h) / h$. Furthermore, we have from the geometrical construction and the assumed boundary resolution that only a few elements $K\in \ThSp$ are mapped to the same element $K'\in\ThS\setminus\ThSp$. 
\end{remark}

The key result for the ghost-penalty operator is the following ghost-penalty mechanism.
\begin{lemma}\label{lem:gp-mechanism}
  Under \Cref{assumption:ghost-penalty} and with the ghost-penalty parameter choice \eqref{eqn.gp-parameter-choice}, it holds for $u_h\in V_h^n$ that
  \begin{equation}\label{eqn:gp-mechanism}
    \clem[a]{lem:gp-mechanism}\left(\norm{\nabla u_h}{\Ohn}^2 + s_h^n(u_h, u_h)\right) \leq \norm{\nabla u_h}{\OdhT}^2 \leq \clem[b]{lem:gp-mechanism}\left(\norm{\nabla u_h}{\Ohn}^2 + s_h^n(u_h, u_h)\right),
  \end{equation}
  with constants $\clem[a]{lem:gp-mechanism}, \clem[b]{lem:gp-mechanism}> 0$ independent of the mesh size and mesh-interface cut configuration.
\end{lemma}
See \cite[Lemma 5.5]{LO19} for a proof thereof.

\begin{remark}
Note that we can choose a smaller, necessary, set of facets for the ghost-penalty stabilisation form, to avoid couplings between all elements with a facet in $\Fhn$, considered for example in~\cite{BNV22, HLSvW22}. However, we keep the above, sufficient, set $\Fhn$ here for simplicity.
\end{remark}

\subsection{Geometry description}
\label{sec:discrete.subsec:geo}
We assume that the discrete geometry $\Ohn$ is represented by a discrete level set function. That is,
\begin{equation*}
  \Ohn \coloneqq \{\xb\in\RR^d:\phi_h(\xb, t^n) \leq 0\},
\end{equation*}
where $\phi_h$ is a piecewise-linear interpolation of a smooth level set function $\phi$. See also \Cite[Section~4]{BCH+14} for further details of geometry descriptions in CutFEM. The piecewise linear level set approximation introduces a geometry approximation error of order $h^2$. However, we will only consider piecewise-linear finite elements, i.e., $k=1$; consequently, second-order convergence is optimal. Higher-order geometry approximations for level set geometries can be achieved in a number of different methods, see, e.g.,~\cite{Leh16, FO15, MKO13, OS16, Say15}. In particular, the isoparametric approach~\cite{Leh16} has also been analysed in the Eulerian moving domain setting~\cite{LL21}. We note that because our method evaluates every $\uhn$ only on $\Ohn$, we expect the isoparametric geometry approximation analysis for this method to be significantly simpler than in~\cite{LL21}, where functions defined with respect to a mesh deformation have to be evaluated accurately on a mesh with a different deformation.

\subsection{Stability}
\label{sec:stability}
Using the ghost-penalty mechanism, we obtain a similar stability result as for the temporally semi-discrete scheme.
\begin{lemma}[Discrete stability for BDF1]\label{lem:stability.discrete.bdf1}
For the solution of the finite element method \eqref{eqn:discrete-scheme.bdf1}, we have for sufficiently small $\dt$ and $h^2$, the stability estimate
\begin{equation*}
  \norm{u_h^N}{\O_h^N}^2 + \dt\sum_{n=1}^N\frac{\nu}{2\clem[b]{lem:gp-mechanism}}\norm{\nabla\uhn}{\OdhT}^2 \leq \exp(\clem{lem:stability.discrete.bdf1}T)\left(\norm{u_h^{0}}{\O_h^{0}}^2 + \dt\sum_{n=1}^{N}\frac{2}{\nu}\norm{f_h^n}{H^{-1}(\Ohn)}^2 \right),
\end{equation*}
with a constant $\clem{lem:stability.discrete.bdf1}$ independent of $h$, $\dt$, $N$ and the mesh-interface cut configurations.
\end{lemma}

\begin{proof}
Testing \eqref{eqn:discrete-scheme.bdf1} with $\vh = 2\dt\uhn$ gives
\begin{equation*}
  2\norm{\uhn}{\Ohn}^2 - 2(\uh^{n-1}, \uhn)_{\Oh^{n-1}} + 2\nu\dt\norm{\nabla\uhn}{\Ohn}^2 - 2\dt(\uhn,\wb^n\cdot\nabla\uhn)_{\Ohn} + 2\dt s_h^n(\uhn,\uhn) = 2\dt(f_h^n,\uhn)_{\Ohn}.
\end{equation*}
As in the proof of \Cref{lemma.smooth-euler.stability}, we use the Cauchy-Schwarz inequality to get
\begin{multline*}
  \norm{\uhn}{\Ohn}^2 + \norm{\uhn - \uh^{n-1}}{\Oh^{n-1}}^2 + 2\nu\dt\norm{\nabla\uhn}{\Ohn}^2 + 2\dt s_h^n(\uhn,\uhn) \\
  \leq \norm{\uh^{n-1}}{\Oh^{n-1}}^2 + 2\dt\norm{f_h^n}{H^{-1}(\Ohn)}\norm{\nabla\uhn}{\Ohn} + 2\dt\sqrt{2} w_\infty\norm{\uhn}{\Ohn}\norm{\nabla\uhn}{\Ohn} + \norm{\uhn}{\Oh^{n-1}\setminus\Ohn}.
\end{multline*}
Dropping the positive term $\norm{\uhn - \uh^{n-1}}{\Oh^{n-1}}^2$ on the left-hand side, using two weighted Young's inequalities and taking the $\norm{\nabla\uhn}{\Ohn}^2$ terms to the left-hand side, as in the proof of \Cref{lemma.smooth-euler.stability}, and using the ghost-penalty mechanism \eqref{eqn:gp-mechanism}, we get
\begin{equation}\label{eqn:stabilit-proof1}
  \norm{\uhn}{\Ohn}^2 + \frac{\nu\dt}{\clem[b]{lem:gp-mechanism}}\norm{\nabla\uhn}{\OdhT}^2 \leq \norm{\uh^{n-1}}{\Oh^{n-1}}^2 + \frac{2\dt}{\nu}\norm{f_h^n}{H^{-1}(\Ohn)}^2 + \frac{4\dt}{\nu}w_\infty^2\norm{\uhn}{\Ohn}^2 + \norm{\uhn}{\Oh^{n-1}\setminus\Ohn}^2
\end{equation}
To bound the solution on the extension strip, we use a discrete version of \eqref{eqn:smooth-extensio-strip-bound}. In particular, it holds under \Cref{assumption:ghost-penalty}, that for all $u_h\in V_h^n$ and any $\eps>0$
\begin{equation}\label{eqn:strip-estimate-discrete}
  c_1 \norm{u_h}{\Sdhp(\Ohn)}^2 \leq \delta_h(1 + \eps^{-1})\norm{u_h}{\Ohn}^2 + \delta_h\eps\norm{\nabla u_h}{\Ohn}^2 + \delta_h L((1 + \eps^{-1})h^2 + \eps)s_h^n(u_h, u_h),
\end{equation}
with $c_1>0$ independent of $h$ and mesh-interface cut configuration, see \cite[Lemma~5.7]{LO19}. With the choice of $\eps < \nu c_1 (2c_{\delta_h}w_{\nb}^\infty\clem[b]{lem:gp-mechanism})^{-1}$ and $h^2$ sufficiently small, we have
\begin{equation*}
  \norm{\uhn}{\Oh^{n-1}\setminus\Ohn}^2
  \leq \norm{\uhn}{\Sdhp(\Ohn)}^2
  \leq c_2 \dt \norm{\uhn}{\Ohn} +  \frac{\nu\dt}{2\clem[b]{lem:gp-mechanism}}\norm{\nabla\uhn}{\OdhT}^2,
\end{equation*}
with $c_2>0$ independent of $h$ and mesh-interface cut configuration. Inserting this in \eqref{eqn:stabilit-proof1} then gives
\begin{equation}\label{eqn:stabilit-proof2}
  \norm{\uhn}{\Ohn}^2 + \frac{\nu\dt}{2\clem[b]{lem:gp-mechanism}}\norm{\nabla\uhn}{\OdhT}^2 \leq \norm{\uh^{n-1}}{\Oh^{n-1}}^2 + \frac{2\dt}{\nu}\norm{f_h^n}{H^{-1}(\Ohn)}^2 + \dt\left(c_2 + \frac{4}{\nu}w_\infty^2\right)\norm{\uhn}{\Ohn}^2.
\end{equation}
Summing over $n=1,\dots, N$ and applying a discrete Grönwall's estimate for $\dt \lesssim \nu w_\infty^{-2}$ sufficiently small gives the result.
\end{proof}

We can also prove a similar result for the BDF2 method \cref{eqn.discrete-scheme:bdf2}. First, we define the discrete BDF2 tuple-norm as
\begin{equation}\label{eqn.def-tuple-norm}
  \norm{(\uhn, \uh^{n-1})}{n}^2\coloneqq
  \norm{\uhn}{\Ohn}^2 + \norm{2\uhn - \uh^{n-1}}{\Oh^{n-1}}^2.
\end{equation}
This definition follows the spirit of only evaluating (discrete) functions on their original domains and domains from previous time steps, rather than domains from future time steps. As this norm and our time-stepping method uses domains from different time step. The polarisation identity leads to the inequality from the following lemma.
\begin{lemma}\label{lem:bdf2-polar}
Let $\uh\in V_h^n$, $\vh\in V_h^{n-1}$ and $\wh\in V_h^{n-2}$. Then, with the BDF2 tuple-norm, we have for any $\eps_1, \eps_2>0$ and constants $\clem[a]{lem:bdf2-polar}, \clem[b]{lem:bdf2-polar}>0$ independent of $\dt$, $h$ and the mesh-interface cut configuration
\begin{multline}\label{eqn:bdf2-polarisation-est}
  \frac{3}{2}(\uh, 4\uh)_{\Ohn} - 2(\vh, 4\uh)_{\Oh^{n-1}} + \frac{1}{2}(\wh, 4\uh)_{\Oh^{n-2}}\\
  \geq
  \begin{aligned}[t]
    &\norm{(\uh, \vh)}{n}^2 - \norm{(\vh, \wh)}{n-1}^2\\
    &- \clem[a]{lem:bdf2-polar}\delta_h (1 + \eps_1^{-1})\norm{\uh}{\Ohn}^2
     - \clem[a]{lem:bdf2-polar}\delta_h \eps_1 \norm{\nabla\uh}{\Ohn}^2
     - \clem[a]{lem:bdf2-polar}\delta_h K ((1 + \eps_1^{-1})h^2 + \eps_1)s_h^n(\uh, \uh)\\
    &- \clem[b]{lem:bdf2-polar}\delta_h (1 + \eps_2^{-1})\norm{\vh}{\Oh^{n-1}}^2
     - \clem[b]{lem:bdf2-polar}\delta_h \eps_2 \norm{\nabla\vh}{\Oh^{n-1}}^2
     - \clem[b]{lem:bdf2-polar}\delta_h K ((1 + \eps_2^{-1})h^2 + \eps_2)s_h^{n-1}(\vh, \vh).
  \end{aligned}
\end{multline}
\end{lemma}

\begin{proof}
We first observe
\begin{align*}
  \norm{\uh - 2\vh + \wh}{\Oh^{n-2}}^2 
    &= \norm{\uh}{\Oh^{n-2}}^2 - 4(\uh,\vh)_{\Oh^{n-2}} + 4\norm{\vh}{\Oh^{n-2}}^2 + 2(\uh,\wh)_{\Oh^{n-2}} - 4(\vh, \wh)_{\Oh^{n-2}} + \norm{\wh}{\Oh^{n-2}}^2,\\
  \intertext{and}
  4(\uh, \vh)_{S} &= 2\norm{\uh}{S}^2 + 2\norm{\vh}{S}^2 - 2\norm{\uh - \vh}{S}^2.
\end{align*}
Consequently, we may rearrange terms to find
\begin{align*}
  6(\uh, \uh)_{\Ohn} &-8(\vh,\uh)_{\Oh^{n-1}} + 2(\wh, \uh)_{\Oh^{n-2}}\\
    &=  \begin{multlined}[t]
            \norm{\uh}{\Ohn}^2 + \norm{2\uh - \vh}{\Oh^{n-1}}^2
          + \norm{\vh}{\Oh^{n-1}}^2 - \norm{2\vh - \wh}{\Oh^{n-2}}^2
          + \norm{\uh - 2\vh + \wh}{\Oh^{n-2}}^2\\
          + 4\Big(\norm{\uh}{\Ohn}^2 - \norm{\uh}{\Oh^{n-1}}^2\Big)
          +  \Big(\norm{\uh}{\Ohn}^2 - \norm{\uh}{\Oh^{n-2}}^2\Big)
          + 4\Big(  (\uh, \vh)_{\Oh^{n-2}} - (\uh, \vh)_{\Oh^{n-1}}\Big)
        \end{multlined}\\
    &=  \begin{multlined}[t]
            \norm{\uh}{\Ohn}^2 + \norm{2\uh - \vh}{\Oh^{n-1}}^2
          - \norm{\vh}{\Oh^{n-1}}^2 - \norm{2\vh - \wh}{\Oh^{n-2}}^2
          + \norm{\uh - 2\vh + \wh}{\Oh^{n-2}}^2 + 2\norm{v_h}{\Oh^{n-2}}\\
          + 5\Big(\norm{\uh}{\Ohn}^2 - \norm{\uh}{\Oh^{n-1}}^2\Big)
          +  \Big(\norm{\uh}{\Oh^{n-2}}^2 - \norm{\uh}{\Oh^{n-1}}^2\Big)
          + 2\Big(\norm{\uh - \vh}{\Oh^{n-1}}^2 - \norm{\uh - \vh}{\Oh^{n-2}}^2\Big).
        \end{multlined}
\end{align*}
Now the first four terms correspond to the desired tuple norm \eqref{eqn.def-tuple-norm}. The next two terms are positive and can be dropped. Furthermore, observing that $\norm{\uh}{S}^2 - \norm{\uh}{K}^2 = \norm{\uh}{S\setminus K}^2 - \norm{\uh}{K\setminus S}^2$, we find with a triangle and Young's inequality for the final term and after dropping positive terms on the left-hand side that
\begin{align}\label{eqn:proof.polar-ineq.1}
  6(\uh, \uh)_{\Ohn}& - 8(\vh,\uh)_{\Oh^{n-1}} + 2(\wh, \uh)_{\Oh^{n-2}}\nonumber\\
      &\geq
      \begin{multlined}[t]
         \norm{(\uh, \vh)}{n}^2 - \norm{(\vh, \wh)}{n-1}^2
        + 5\norm{\uh}{\Ohn\setminus\Oh^{n-1}}^2
        - 5\norm{\uh}{\Oh^{n-1}\setminus\Ohn}^2\\
        +  \norm{\uh}{\Oh^{n-2}\setminus\Oh^{n-1}}^2
        -  \norm{\uh}{\Oh^{n-1}\setminus\Oh^{n-2}}^2
        + 2\norm{\uh - \vh}{\Oh^{n-1}\setminus\Oh^{n-2}}^2 
        - 2\norm{\uh - \vh}{\Oh^{n-2}\setminus\Oh^{n-1}}^2
      \end{multlined}\nonumber\\
      &\geq
        \norm{(\uh, \vh)}{n}^2 - \norm{(\vh, \wh)}{n-1}^2
        - 5\norm{\uh}{\Oh^{n-1}\setminus\Ohn}^2
        -  \norm{\uh}{\Oh^{n-1}\setminus\Oh^{n-2}}^2
        - 4\norm{\uh}{\Oh^{n-2}\setminus\Oh^{n-1}}^2
        - 4\norm{\vh}{\Oh^{n-2}\setminus\Oh^{n-1}}^2.
\end{align}
The final four norms in \eqref{eqn:proof.polar-ineq.1} are over subsets of the strips $\Sdh$ around $\Ohn$ and $\Oh^{n-1}$, respectively. Therefore, the result follows by applying \eqref{eqn:strip-estimate-discrete}.
\end{proof}

\begin{lemma}[Discrete stability for BDF2]\label{lem:stability.discrete.bdf2}
Let the fully discrete BDF2 scheme \eqref{eqn.discrete-scheme:bdf2} be initialised by a single step of the fully discrete BDF1 method \eqref{eqn:discrete-scheme.bdf1}. For $\dt$ and $h^2$ sufficiently small, we then have the stability estimate for the BDF2 scheme
\begin{equation*}
  \norm{(\uh^N, \uh^{N-1})}{N}^2 
  + \dt\sum_{n=1}^{N}\frac{\nu}{2\clem[b]{lem:gp-mechanism}}\norm{\nabla\uhn}{\OdhT}^2
    \leq \exp(\clem{lem:stability.discrete.bdf2}T)\left[6\norm{\uh^{0}}{\Oh^0} + \dt\sum_{n=1}^N\frac{10}{\nu} \norm{f_h^n}{H^{-1}(\Ohn)}^2\right],
\end{equation*}
with a constant $\clem{lem:stability.discrete.bdf2}$ independent of $h$, $\dt$, $N$ and the mesh-interface cut configurations.
\end{lemma}

\begin{proof}
Testing \eqref{eqn.discrete-scheme:bdf2} with $\vh=4\dt\uhn$, using the Cauchy-Schwarz inequality, weighted Young's inequality and \eqref{eqn:bdf2-polarisation-est}, we get
\begin{multline*}
  \norm{(\uhn, \uh^{n-1})}{n}^2 
  + \dt2(\nu\norm{\nabla\uhn}{\Ohn}^2 + \nu s_h^n(\uhn,\uhn))
  \leq \norm{(\uh^{n-1}, \uh^{n-2})}{n-1}^2
    + \frac{8\dt}{\nu}w_\infty^2\norm{\uhn}{\Ohn}^2
    + \frac{4\dt}{\nu} \norm{f_h^n}{H^{-1}(\Ohn)}^2\\
    + c\dt (1 + \eps_1^{-1})\norm{\uhn}{\Ohn}^2
    + c\dt \eps_1 \norm{\nabla\uhn}{\Ohn}^2
    + c\dt L ((1 + \eps_1^{-1})h^2 + \eps_1)s_h^n(\uhn, \uhn)\\
    + c\dt (1 + \eps_2^{-1})\norm{\uh^{n-1}}{\Oh^{n-1}}^2
    + c\dt \eps_2 \norm{\nabla\uh^{n-1}}{\Oh^{n-1}}^2
    + c\dt L ((1 + \eps_2^{-1})h^2 + \eps_2)s_h^{n-1}(\uh^{n-1}, \uh^{n-1}).
\end{multline*}
Using the ghost-penalty mechanism and choosing $\eps_1,\eps_2$ sufficiently small, we have for $h^2$ sufficiently small
\begin{multline*}
  \norm{(\uh, \uh^{n-1})}{n}^2 
  + \dt\frac{\nu}{\clem[b]{lem:gp-mechanism}}\norm{\nabla\uhn}{\OdhT}^2
    \leq \norm{(\uh^{n-1}, \uh^{n-2})}{n-1}^2
    + \dt(c_1(\nu^{-1}) + \frac{8}{\nu} w_\infty^2)\norm{\uhn}{\Ohn}^2
    + \frac{4\dt}{\nu} \norm{f_h^n}{H^{-1}(\Ohn)}^2\\
    + c_2(\nu^{-1})\dt \norm{\uh^{n-1}}{\Oh^{n-1}}^2
    + \dt\frac{\nu}{2\clem[b]{lem:gp-mechanism}}\norm{\nabla\uh^{n-1}}{\OdhTlast}^2.
\end{multline*}
Summing over $n=2,\dots,N$, and using the triangle inequality for the BDF2 tuple norm
\begin{multline*}
  \norm{(\uh^N, \uh^{N-1})}{N}^2 
  + \dt\sum_{n=2}^{N}\frac{\nu}{2\clem[b]{lem:gp-mechanism}}\norm{\nabla\uhn}{\OdhT}^2
    \leq 5\norm{\uh^{1}}{\Oh^1} + \norm{\uh^{0}}{\Oh^2}^2
    + \frac{\dt\nu}{\clem[b]{lem:gp-mechanism}}\norm{\nabla\uh^{1}}{\OdhTone}^2\\
    + \dt\sum_{n=1}^N\big(c_1(\nu^{-1}) + c_2(\nu^{-1}) + \frac{8}{\nu} w_\infty^2\big)\norm{\uhn}{\Ohn}^2
    + \dt\sum_{n=2}^N\frac{4}{\nu} \norm{f_h^n}{H^{-1}(\Ohn)}^2.   
\end{multline*}
We then add $\frac{\dt\nu}{2\clem[b]{lem:gp-mechanism}}\norm{\nabla\uh^{1}}{\OdhTone}^2$ to both sides and use that $\uh^1$ is the solution from a single step of the BDF1 method, which allows us to bound the $\uh^1$ terms by the initial condition and data using
\eqref{eqn:stabilit-proof2} with $n=1$. Finally, observing that $\norm{\uhn}{\Ohn}^2 \leq \norm{(\uhn,\uh^{n-1})}{n}^2$, the result follows for $\dt$ sufficiently small from a discrete Grönwall inequality.
\end{proof}

\begin{remark}[Error analysis]
Following standard lines of argument for the consistency error of our method requires a suitable bound of
\begin{equation*}
  \dod[2]{}{t}\int_{\O(t)}u(t)\vh\dif x,
\end{equation*}
with $v_h$ from the finite element space.
Using Reynolds transport theorem, cf.\ \eqref{eqn.identity1} as well as the hyper-surface version thereof, see, e.g., \cite[Lemma~2.1]{DE12} leads to the identity
\begin{align*}
  \dod[2]{}{t}\int_{\O(t)}u(t)\vh\dif x
  &= \dod{}{t}\left(\int_{\O(t)}\vh\partial_t u\dif x 
    + \int_{\G(t)} \wb\cdot\bn (u \vh) \dif s\right) \\
  &= \begin{multlined}[t]
      \int_{\O(t)}\vh\partial_{tt} u\dif x
      +\int_{\G(t)} \vh \wb\cdot\bn \partial_t u \dif s\\
      \int_{\G(t)} \vh \partial_t(\wb\cdot\bn u) 
        + \wb\cdot\nabla (\wb\cdot\bn u \vh)
        + \vh \wb\cdot\bn u\div_\Gamma\wb \dif s.
    \end{multlined}
\end{align*}
To bound the right-hand side of this, we require control of $\wb\cdot\nabla (\wb\cdot\bn u \vh)$ on $\G(t)$, i.e., we need $v_h\in H^{3/2}$, which is not given. Consequently, deriving an error estimate for our suggested scheme is not obvious, even for the BDF1 case, and remains on open problem.
\end{remark}

\section{Numerical Examples}
\label{sec:num-ex}

We implement the method using \texttt{ngsxfem}~\cite{LHPvW21}, an add-on to \texttt{Netgen/NGSolve}~\cite{Sch97, Sch14} for unfitted finite element discretisations. Throughout, for a given spatial norm $\Vert\cdot\Vert_X$, we measure the error in the discrete space-time norm
\begin{equation}
  \Vert \cdot \Vert^2_{ L^2(X)} = \sum_{i=1}^n \dt \Vert \cdot\Vert_X^2.
\end{equation}
Throughout our numerical examples, we shall consider $k=1$.

\subsection{Example 1: Travelling circle}
\label{sec:num-ex.subsec:travelling}

As our first example, we consider the simple geometry of a travelling circle, taken from~\cite{LO19}.

\subsubsection*{Set-up}
We consider the background domain $\widetilde{\O} = (-0.7, 0.9)\times(-0.7, 0.7)$ over the time interval $[0, 0.2]$. The geometry and transport-field are given through
\begin{equation}
  \rho(\xb, t) = \big((\sin(2\pi t) / \pi), 0\big)^T,\quad
  \phi(\xb, t) = \Vert \xb - \rho(\xb, t) \Vert_2 - 0.5,\quad
  \wb(\xb, t) = \partial_t \rho(\xb, t),
\end{equation}
and the exact solution is set as
\begin{equation}
  u(\xb, t) = \cos^2\big(\pi \Vert \xb - \rho(\xb, t)\Vert_2\big).
\end{equation}
The viscosity is chosen as $\nu=1$, and the forcing term is set according to \eqref{eqn.strong-problem.1}. The initial mesh size is $h_0=0.4$, and the initial time step is $\dt_0=0.1$.

\subsubsection*{Convergence Study}

We consider a series of uniform mesh and time-step refinements. To this end, we construct a series of meshes with $h= 2^{-L_x}h_0$ together with time steps $\dt = 2^{-L_t}\dt_0 $, for $L_x=0,\dots, 5$ and $L_t=0,\dots, 6$. Experimental order of convergence rates in the mesh size $\eoc{x}$ and in time step $\eoc{t}$ are then computed by comparing the resulting errors between two levels, where the other parameter (time step or mesh size) is the most refined available. Combined rates ($\eoc{xt}$) are computed by comparing vales after refining both the mesh size and time step.

The resulting errors of the BDF1 scheme can be seen in \Cref{tab:travelling-cir.bdf1.l2l2} for the $ L^2(L^2)$-norm, in \Cref{tab:travelling-cir.bdf1.linfl2} for the $L^\infty(L^2)$-norm and in \Cref{tab:travelling-cir.bdf1.l2h1} for the $ L^2(H^1)$-norm. In the $ L^2(L^2)$-norm, we observe optimal linear convergence in time on the finest mesh ($\eoc{t}\approx 1$) and close to second-order convergence in space with the smallest time step ($\eoc{x}\approx 2$). The errors in the $L^\infty(L^2)$-norm are very similar to those in the $L^2(L^2)$-norm, with linear convergence in time on the finest mesh and initial second-order convergence in space using the smallest time step. In the $ L^2(H^1)$-norm, we observe no further time step convergence after the first refinement on the finest mesh. However, we observe optimal, linear convergence for both mesh refinement with the smallest time step ($\eoc{x}\approx 1$) and for combined mesh and time-step refinement ($\eoc{xt}\approx 1$).

The errors resulting from the BDF2 scheme in the $ L^2(L^2)$-, $L^\infty(L^2)$- and $ L^2(H^2)$-norm are presented in \Cref{tab:travelling-cir.bdf2.l2l2}, \Cref{tab:travelling-cir.bdf2.linfl2} and \Cref{tab:travelling-cir.bdf2.l2h1}, respectively. In the $L^2(L^2)$-norm, we observe between first and second-order convergence with respect to the time step on the fines mesh. Nevertheless, we observe optimal second-order convergence with respect to the mesh size with the finest mesh size and for combined mesh and time-step refinement. As in the BDF1 case, the $L^\infty(L^2)$-norm results are very similar to the $L^2(L^2)$ results, with some suboptimal convergence in time on the finest mesh, but optimal rates in space and for combined refinement ($\eoc{x}\approx\eoc{xt}\approx2$). Concerning the error in the $ L^2(H^1)$-norm, it appears that the spatial error again dominated with very little time step convergence on the finest mesh. Nevertheless, we again observe optimal linear convergence in space with the smallest time step and optimal second-order convergence when refining the time step once and the mesh twice ($\eoc{xxt}\approx 2$), i.e., comparing the error from $(L_t, L_x)$ with that from $(L_t - 1, L_x-2)$.

\begin{figure}
  \centering
  \includegraphics{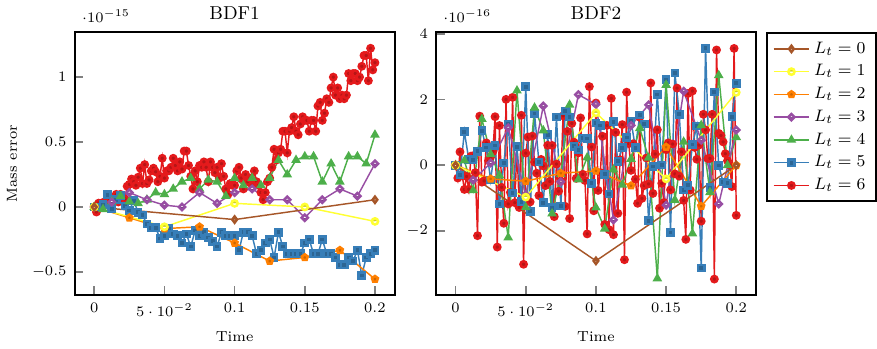}
  \caption{Mass conservation error for the travelling circle example on the coarsest mesh $L_x=0$.}
  \label{fig.mass-conservation}
\end{figure}

Finally, we look at the conservation of mass resulting from the scheme in \Cref{fig.mass-conservation} for both the BDF1 and BDF2 schemes on the coarsest mesh. We see that mass is conserved up to machine precision for both schemes, as expected.

\begin{table}
\centering
\caption{Convergence results for the BDF1 scheme in the $ L^2(L^2)$ norm for the translating circle test case in \Cref{sec:num-ex.subsec:travelling}}
\begin{tabular}{llllllll}
\toprule
$L_t$~\textbackslash~$L_x$  & 0 & 1 & 2 & 3 & 4 & 5 & $\eoc{t}$ \\
\midrule
0 & \num{1.17e-01} & \num{4.10e-02} & \num{2.04e-02} & \num{1.54e-02} & \num{1.39e-02} & \num{1.32e-02} & --\\
1 & \uline{\num{1.06e-01}} & \num{3.57e-02} & \num{1.34e-02} & \num{8.42e-03} & \num{7.18e-03} & \num{6.79e-03} & 0.96\\
2 & \num{1.00e-01} & \uline{\num{3.32e-02}} & \num{1.09e-02} & \num{4.88e-03} & \num{3.74e-03} & \num{3.47e-03} & 0.97\\
3 & \num{9.72e-02} & \num{3.18e-02} & \uline{\num{9.64e-03}} & \num{3.38e-03} & \num{2.01e-03} & \num{1.78e-03} & 0.96\\
4 & \num{9.57e-02} & \num{3.12e-02} & \num{9.03e-03} & \uline{\num{2.70e-03}} & \num{1.20e-03} & \num{9.27e-04} & 0.94\\
5 & \num{9.50e-02} & \num{3.09e-02} & \num{8.74e-03} & \num{2.39e-03} & \uline{\num{8.20e-04}} & \num{5.03e-04} & 0.88\\
6 & \num{9.46e-02} & \num{3.07e-02} & \num{8.60e-03} & \num{2.26e-03} & \num{6.55e-04} & \uline{\num{2.96e-04}} & 0.76\\
\midrule
$\eoc{x}$ & -- & 1.62 & 1.84 & 1.93 & 1.79 & 1.14 \\
$\uline{\eoc{xt}}$ & -- & 1.68 & 1.78 & 1.84 & 1.72 & 1.47 \\
\bottomrule
\end{tabular}
\label{tab:travelling-cir.bdf1.l2l2}
\end{table}

\begin{table}
\centering
\caption{Convergence results for the BDF1 scheme in the $ L^\infty(L^2)$ norm for the translating circle test case in \Cref{sec:num-ex.subsec:travelling}}
\begin{tabular}{llllllll}
\toprule
$L_t$\textbackslash $L_x$  & 0 & 1 & 2 & 3 & 4 & 5 & $\text{eoc}_{t}$ \\
\midrule
0 & \num{3.15e-01} & \num{1.04e-01} & \num{4.73e-02} & \num{3.58e-02} & \num{3.27e-02} & \num{3.16e-02} & --\\
1 & \uline{\num{3.17e-01}} & \num{9.62e-02} & \num{3.26e-02} & \num{1.99e-02} & \num{1.77e-02} & \num{1.71e-02} & 0.88\\
2 & \num{3.15e-01} & \uline{\num{9.43e-02}} & \num{2.76e-02} & \num{1.17e-02} & \num{9.34e-03} & \num{8.99e-03} & 0.93\\
3 & \num{3.15e-01} & \num{9.32e-02} & \uline{\num{2.57e-02}} & \num{8.17e-03} & \num{5.00e-03} & \num{4.64e-03} & 0.95\\
4 & \num{3.15e-01} & \num{9.30e-02} & \num{2.49e-02} & \uline{\num{6.83e-03}} & \num{2.92e-03} & \num{2.39e-03} & 0.96\\
5 & \num{3.15e-01} & \num{9.28e-02} & \num{2.46e-02} & \num{6.34e-03} & \uline{\num{2.01e-03}} & \num{1.26e-03} & 0.92\\
6 & \num{3.15e-01} & \num{9.25e-02} & \num{2.44e-02} & \num{6.17e-03} & \num{1.68e-03} & \uline{\num{7.27e-04}} & 0.80\\
\midrule
$\text{eoc}_{x}$ & -- & 1.77 & 1.92 & 1.98 & 1.88 & 1.21 \\
$\uline{\text{eoc}_{xt}}$ & -- & 1.75 & 1.88 & 1.91 & 1.76 & 1.47 \\
\bottomrule
\end{tabular}
\label{tab:travelling-cir.bdf1.linfl2}
\end{table}

\begin{table}
\centering\caption{Convergence results for the BDF1 scheme in the $ L^2(H^1)$ norm for the translating circle test case in \Cref{sec:num-ex.subsec:travelling}}
\begin{tabular}{llllllll}
\toprule
$L_t$~\textbackslash~$L_x$  & 0 & 1 & 2 & 3 & 4 & 5 & $\eoc{t}$ \\
\midrule
0 & \num{5.60e-01} & \num{3.30e-01} & \num{2.05e-01} & \num{1.34e-01} & \num{9.75e-02} & \num{8.04e-02} & --\\
1 & \uline{\num{5.54e-01}} & \num{3.19e-01} & \num{1.80e-01} & \num{1.06e-01} & \num{6.67e-02} & \num{4.83e-02} & 0.74\\
2 & \num{5.42e-01} & \uline{\num{3.16e-01}} & \num{1.77e-01} & \num{9.33e-02} & \num{5.22e-02} & \num{3.23e-02} & 0.58\\
3 & \num{5.36e-01} & \num{3.15e-01} & \uline{\num{1.76e-01}} & \num{9.13e-02} & \num{4.68e-02} & \num{2.56e-02} & 0.33\\
4 & \num{5.32e-01} & \num{3.14e-01} & \num{1.75e-01} & \uline{\num{9.03e-02}} & \num{4.58e-02} & \num{2.33e-02} & 0.13\\
5 & \num{5.31e-01} & \num{3.14e-01} & \num{1.74e-01} & \num{8.98e-02} & \uline{\num{4.54e-02}} & \num{2.29e-02} & 0.03\\
6 & \num{5.30e-01} & \num{3.14e-01} & \num{1.74e-01} & \num{8.96e-02} & \num{4.52e-02} & \uline{\num{2.27e-02}} & 0.01\\
\midrule
$\eoc{x}$ & -- & 0.76 & 0.85 & 0.96 & 0.99 & 0.99 \\
$\uline{\eoc{xt}}$ & -- & 0.81 & 0.85 & 0.96 & 0.99 & 1.00 \\
\bottomrule
\end{tabular}
\label{tab:travelling-cir.bdf1.l2h1}
\end{table}

\begin{table}
\centering\caption{Convergence results for the BDF2 scheme in the $ L^2(L^2)$ norm for the translating circle test case in \Cref{sec:num-ex.subsec:travelling}}
\begin{tabular}{llllllll}
\toprule
$L_t$~\textbackslash~$L_x$  & 0 & 1 & 2 & 3 & 4 & 5 & $\eoc{t}$ \\
\midrule
0 & \num{1.16e-01} & \num{4.68e-02} & \num{2.57e-02} & \num{1.68e-02} & \num{1.30e-02} & \num{1.13e-02} & --\\
1 & \uline{\num{1.06e-01}} & \num{3.65e-02} & \num{1.41e-02} & \num{7.46e-03} & \num{5.20e-03} & \num{4.35e-03} & 1.38\\
2 & \num{1.01e-01} & \uline{\num{3.35e-02}} & \num{1.02e-02} & \num{3.91e-03} & \num{2.06e-03} & \num{1.51e-03} & 1.53\\
3 & \num{9.78e-02} & \num{3.23e-02} & \uline{\num{9.41e-03}} & \num{2.61e-03} & \num{9.32e-04} & \num{5.06e-04} & 1.58\\
4 & \num{9.60e-02} & \num{3.14e-02} & \num{8.96e-03} & \uline{\num{2.36e-03}} & \num{6.16e-04} & \num{2.12e-04} & 1.26\\
5 & \num{9.51e-02} & \num{3.10e-02} & \num{8.72e-03} & \num{2.25e-03} & \uline{\num{5.67e-04}} & \num{1.45e-04} & 0.54\\
6 & \num{9.47e-02} & \num{3.08e-02} & \num{8.59e-03} & \num{2.19e-03} & \num{5.47e-04} & \uline{\num{1.37e-04}} & 0.09\\
\midrule
$\eoc{x}$ & -- & 1.62 & 1.84 & 1.97 & 2.00 & 2.00 \\
$\uline{\eoc{xt}}$ & -- & 1.66 & 1.83 & 1.99 & 2.06 & 2.05 \\
\bottomrule
\end{tabular}
\label{tab:travelling-cir.bdf2.l2l2}
\end{table}

\begin{table}
\centering
\caption{Convergence results for the BDF2 scheme in the $ L^\infty(L^2)$ norm for the translating circle test case in \Cref{sec:num-ex.subsec:travelling}}
\begin{tabular}{llllllll}
\toprule
$L_t$\textbackslash $L_x$  & 0 & 1 & 2 & 3 & 4 & 5 & $\text{eoc}_{t}$ \\
\midrule
0 & \num{3.11e-01} & \num{1.11e-01} & \num{6.21e-02} & \num{4.36e-02} & \num{3.47e-02} & \num{3.02e-02} & --\\
1 & \uline{\num{3.16e-01}} & \num{9.83e-02} & \num{3.26e-02} & \num{2.13e-02} & \num{1.67e-02} & \num{1.45e-02} & 1.06\\
2 & \num{3.16e-01} & \uline{\num{9.49e-02}} & \num{2.63e-02} & \num{1.04e-02} & \num{6.80e-03} & \num{5.76e-03} & 1.33\\
3 & \num{3.16e-01} & \num{9.44e-02} & \uline{\num{2.54e-02}} & \num{6.57e-03} & \num{2.82e-03} & \num{2.03e-03} & 1.51\\
4 & \num{3.15e-01} & \num{9.32e-02} & \num{2.49e-02} & \uline{\num{6.31e-03}} & \num{1.58e-03} & \num{7.63e-04} & 1.41\\
5 & \num{3.15e-01} & \num{9.31e-02} & \num{2.47e-02} & \num{6.21e-03} & \uline{\num{1.55e-03}} & \num{3.84e-04} & 0.99\\
6 & \num{3.15e-01} & \num{9.29e-02} & \num{2.45e-02} & \num{6.13e-03} & \num{1.53e-03} & \uline{\num{3.81e-04}} & 0.01\\
\midrule
$\text{eoc}_{x}$ & -- & 1.76 & 1.92 & 2.00 & 2.00 & 2.01 \\
$\uline{\text{eoc}_{xt}}$ & -- & 1.74 & 1.90 & 2.01 & 2.03 & 2.02 \\
\bottomrule
\end{tabular}
\label{tab:travelling-cir.bdf2.linfl2}
\end{table}

\begin{table}
\centering
\caption{Convergence results for the BDF2 scheme in the $ L^2(H^1)$ norm for the translating circle test case in \Cref{sec:num-ex.subsec:travelling}}
\begin{tabular}{llllllll}
\toprule
$L_t$~\textbackslash~$L_x$  & 0 & 1 & 2 & 3 & 4 & 5 & $\eoc{t}$ \\
\midrule
0 & \num{5.53e-01} & \num{3.69e-01} & \num{2.48e-01} & \num{1.66e-01} & \num{1.16e-01} & \num{8.87e-02} & --\\
1 & \num{5.49e-01} & \num{3.23e-01} & \num{1.93e-01} & \num{1.16e-01} & \num{7.20e-02} & \num{4.78e-02} & 0.89\\
2 & \dashuline{\num{5.47e-01}} & \num{3.17e-01} & \num{1.77e-01} & \num{9.91e-02} & \num{5.60e-02} & \num{3.24e-02} & 0.56\\
3 & \num{5.39e-01} & \uwave{\num{3.16e-01}} & \dashuline{\num{1.76e-01}} & \num{9.13e-02} & \num{4.86e-02} & \num{2.61e-02} & 0.31\\
4 & \num{5.34e-01} & \num{3.15e-01} & \uline{\num{1.75e-01}} & \uwave{\num{9.07e-02}} & \num{4.57e-02} & \num{2.36e-02} & 0.14\\
5 & \num{5.31e-01} & \num{3.14e-01} & \num{1.75e-01} & \uuline{\num{9.02e-02}} & \uline{\num{4.55e-02}} & \num{2.28e-02} & 0.05\\
6 & \num{5.30e-01} & \num{3.14e-01} & \num{1.74e-01} & \num{8.98e-02} & \num{4.53e-02} & \uuline{\num{2.27e-02}} & 0.00\\
\midrule
$\eoc{x}$ & -- & 0.76 & 0.85 & 0.96 & 0.99 & 1.00 \\
$\uline{\eoc{xxt}}$ & -- & -- & 1.63 & 1.80 & 1.95 & 1.99 \\
\bottomrule
\end{tabular}
\label{tab:travelling-cir.bdf2.l2h1}
\end{table}

\subsection{Example 2: Kite Geometry}
\label{sec:num-ex.subsec:kite}

\begin{figure}
  \centering
  -1 \includegraphics[height=6.2pt, width=0.8\textwidth]{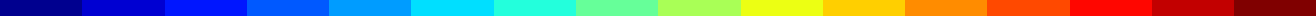} 1\phantom{-}\\[4pt]
  \includegraphics[width=0.19\textwidth]{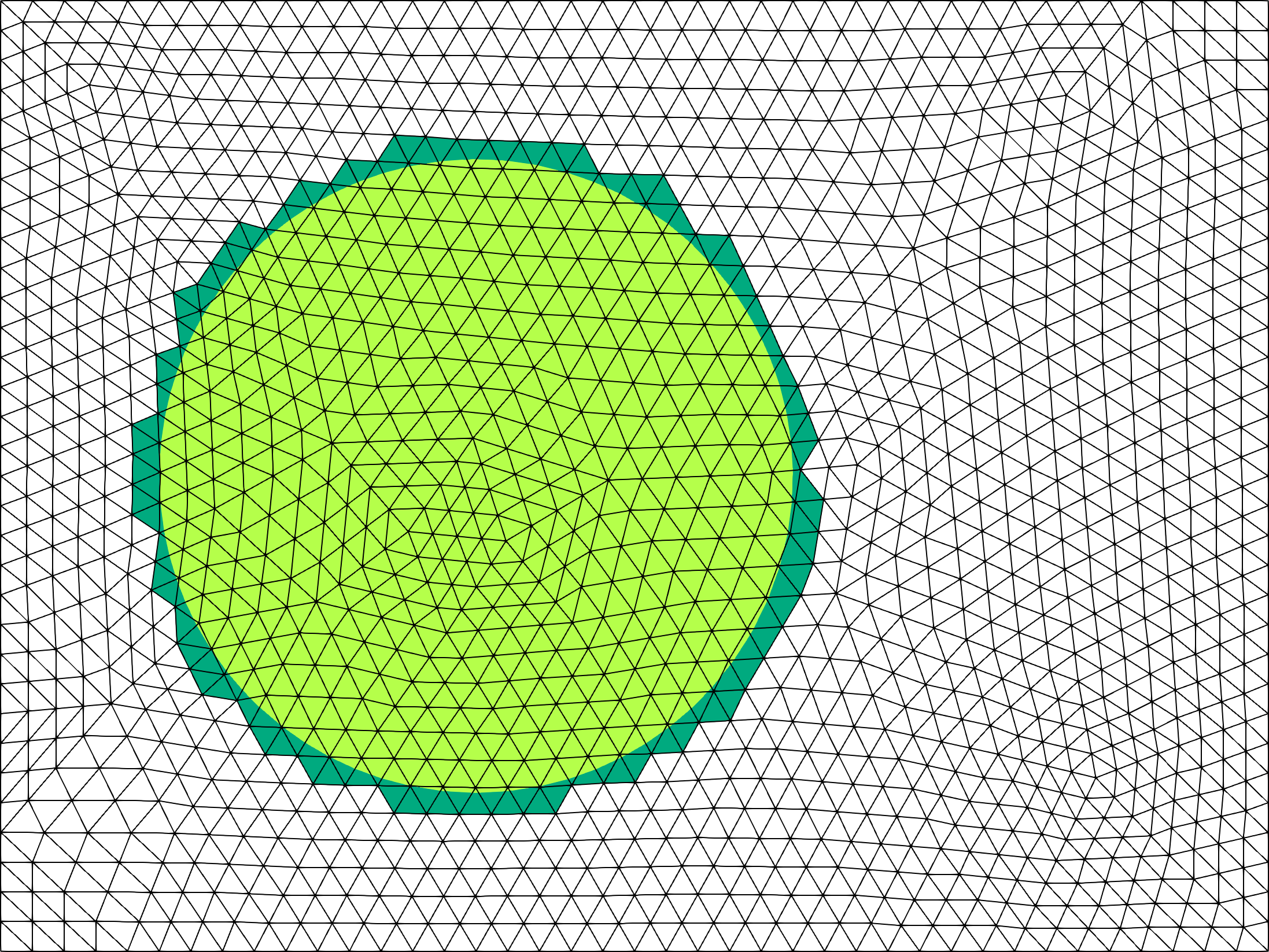}\hspace{1pt}%
  \includegraphics[width=0.19\textwidth]{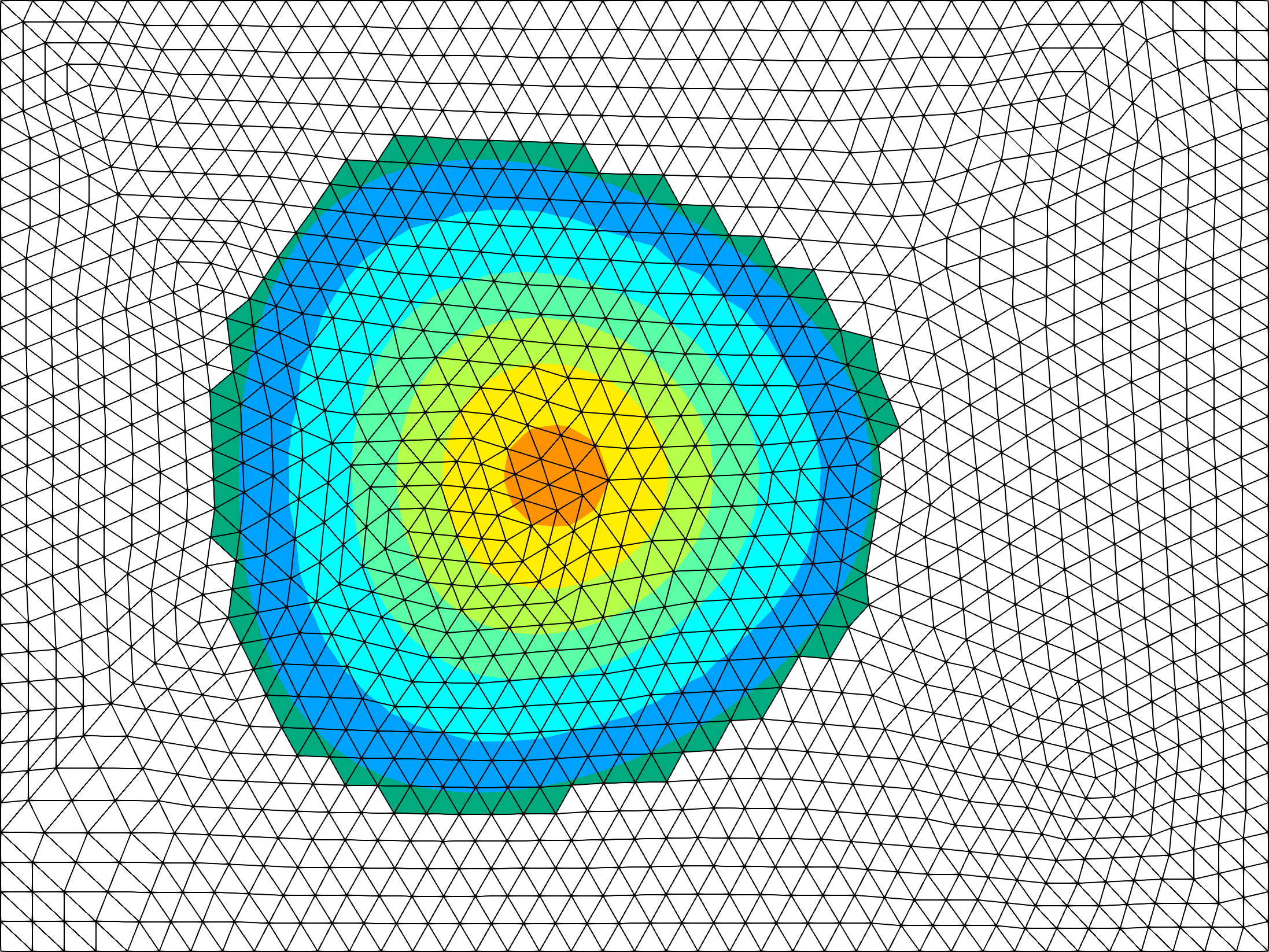}\hspace{1pt}%
  \includegraphics[width=0.19\textwidth]{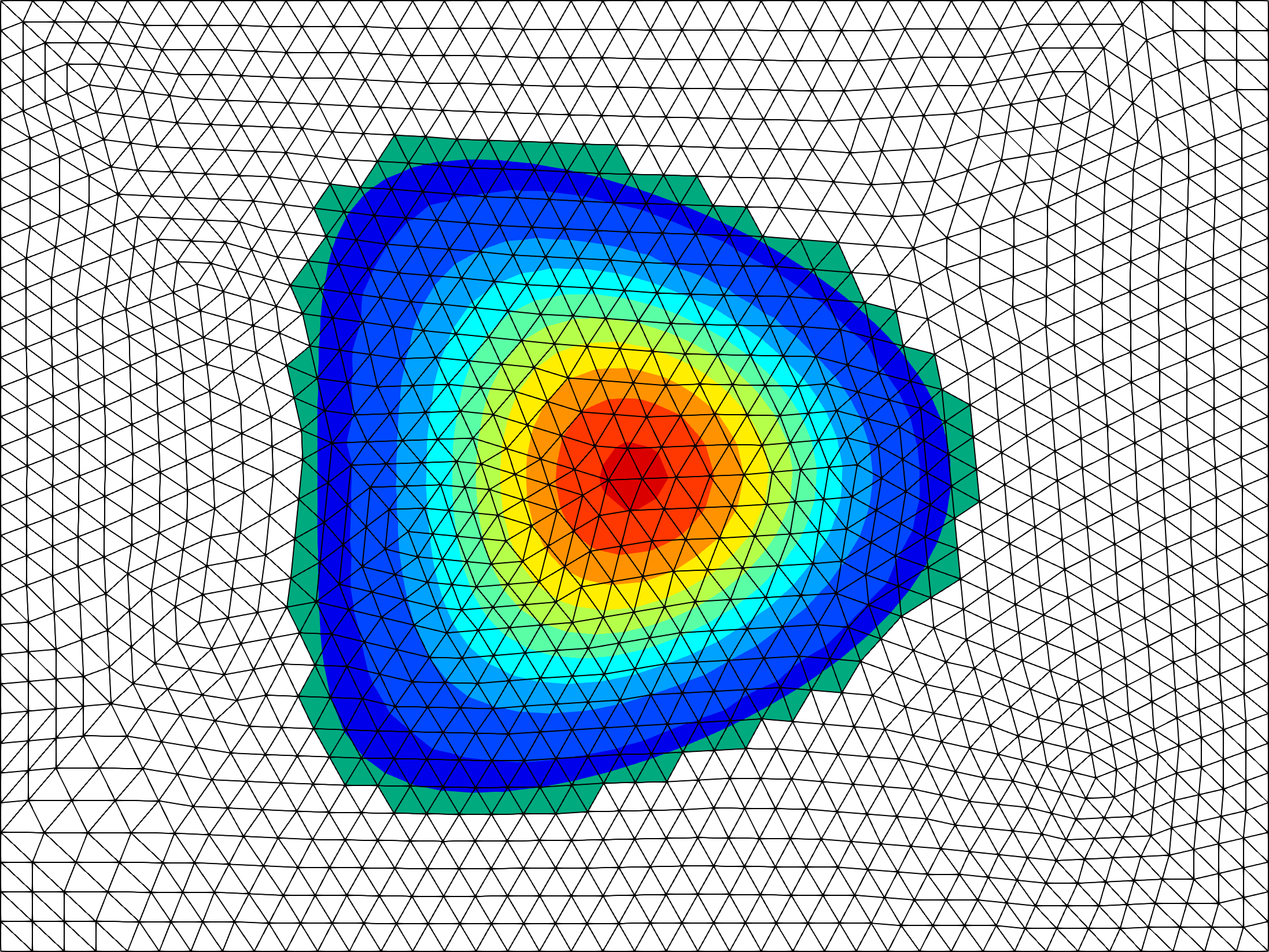}\hspace{1pt}%
  \includegraphics[width=0.19\textwidth]{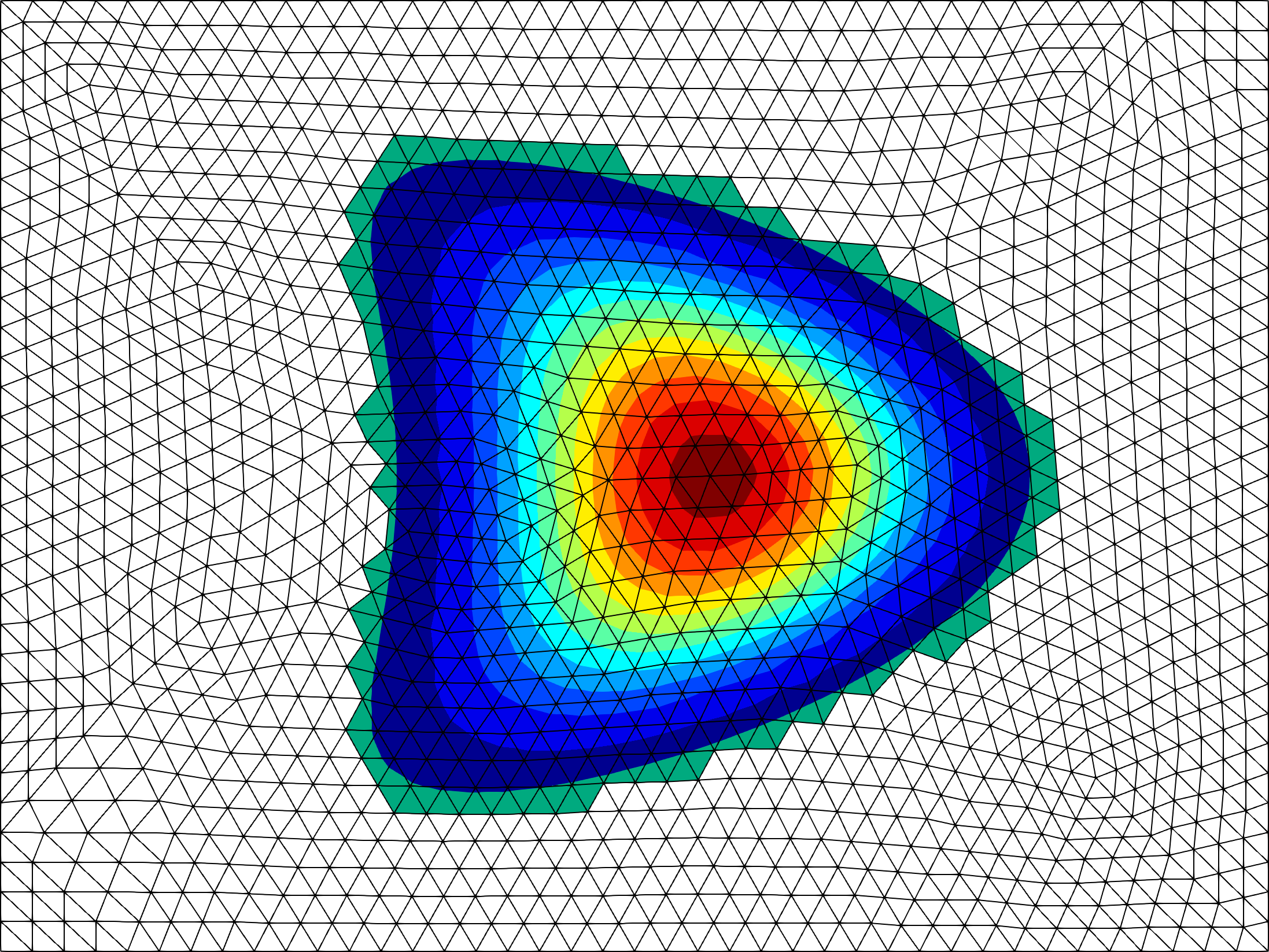}\hspace{1pt}%
  \includegraphics[width=0.19\textwidth]{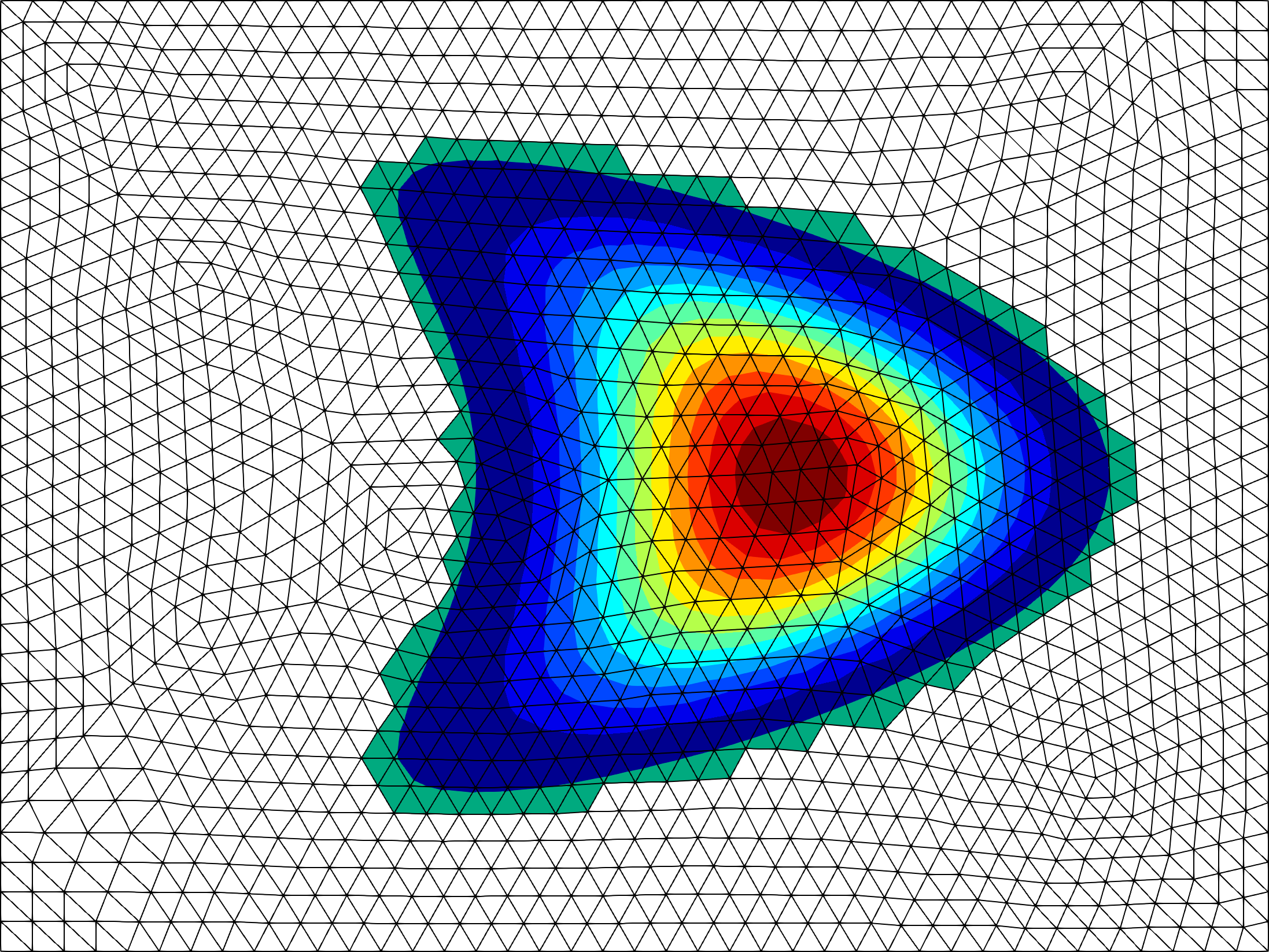}
  \caption{Discrete solution for the kite geometry at times $0.25\tend$ intervals with $h=0.4 \cdot 2^{-2}$ and $\dt=0.5 \cdot 2^{-6}$. Extension elements in $\ThS$ are marked in green.}
  \label{fig:kite}
\end{figure}

As a second example, we consider a case where the geometry translates and deforms in time. This example follows from~\cite{HLP23}. The geometry starts as a circle and changes into a kite-like geometry over time. An illustration of this is given in \Cref{fig:kite}.

\subsubsection*{Set-up}
We take the background domain as $\widetilde\O = (-1.5, 2.5)\times(-1.5, 1.5)$ and $\tend=1$.  The level set function and transport field are given by
\begin{equation*}
  \rho(\xb, t) = \big((1 - \xb_2)^2 t, 0\big)^T, \quad \phi(\xb, t) = \Vert\xb - \rho(\xb, t)\Vert_2 - R, \quad\wb(\xb, t)= \partial_t\rho(\xb, t),
\end{equation*}
with $R=1$. The right-hand side forcing term is taken such that the exact solution is
\begin{equation*}
  u_{ex}(\xb, t) = \cos(\pi\Vert\xb - \rho(\xb, t)\Vert_2 / R ) \sin(t \pi / 2).
\end{equation*}
The viscosity is chosen to be $\nu=0.2$. The initial mesh size and time step are $h=0.4$ and $\dt=0.5$, respectively.

\subsubsection*{Convergence Study}
Again, we consider a series of uniform refinements in space and time. The results for the BDF1 scheme can be found in \Cref{tab:kite.bdf1.l2l2}, \Cref{tab:kite.bdf1.linfl2} and \Cref{tab:kite.bdf1.l2h1} for the $ L^2(L^2)$-, $L^\infty(L^2)$- and $ L^2(H^2)$-norms, respectively. The results for the BDF2 scheme are presented in \Cref{tab:kite.bdf2.l2l2}, \Cref{tab:kite.bdf2.linfl2}, and \Cref{tab:kite.bdf2.l2h1} in the same norms.

For the results of the BDF1 scheme in the $L^2(L^2)$-norm, we again observe optimal linear convergence in time on the finest mesh and linear convergence for combined mesh and time-step refinement, indicating that the temporal error is dominating in this case. In the $L^\infty(L^2)$-norm, we see linear convergence in time, with a small drop in the rate for large time steps, and optimal linear convergence for combined mesh and time-step refinement. In the $ L^2(H^1)$-norm, the convergence with respect to the time step is less than optimal, with both $\eoc{t} = \eoc{xt} \approx 0.7$.

For the BDF2 scheme, we see that for the $L^2(L^2)$-norm error, the order of convergence in time ($\eoc{t}$) starts sub-optimal but increases to $1.84$. This is also reflected in $\eoc{xt}$, which is only around $2$ after the first level of refinement. This suggests that the time step really needs to be taken sufficiently small.
The $L^\infty(L^2)$-norm results are again similar to the $L^2(L^2)$ results. Convergence in time on the finest mesh is slightly suboptimal, with a fraction lower rates than in the $L^2(L^2)$-norm. Similarly, we observe optimal second-order convergence for combined refinement after the first refinement level.
Regarding the error $L^2(H^1)$-norm, we again have sub-optimal convergence in time on the fines mesh with $\eoc{t}\approx1$. However, we see that $\eoc{xxt}$ approaches the optimal value of 2 under refinement. This again suggests that the time step is not sufficiently small on the coarser time levels and that optimal convergence in time can only be realised in combination with mesh refinement.

\begin{table}
\centering
\caption{Convergence results for the BDF1 scheme in the $ L^2(L^2)$ norm for the kite geometry test case in \Cref{sec:num-ex.subsec:kite}}
\begin{tabular}{llllllll}
\toprule
$L_t$~\textbackslash~$L_x$  & 0 & 1 & 2 & 3 & 4 & 5 & $\eoc{t}$ \\
\midrule
0 & \num{6.08e-01} & \num{4.60e-01} & \num{4.25e-01} & \num{4.17e-01} & \num{4.15e-01} & \num{4.14e-01} & --\\
1 & \uline{\num{3.45e-01}} & \num{2.55e-01} & \num{2.27e-01} & \num{2.20e-01} & \num{2.18e-01} & \num{2.18e-01} & 0.93\\
2 & \num{2.50e-01} & \uline{\num{1.45e-01}} & \num{1.23e-01} & \num{1.17e-01} & \num{1.16e-01} & \num{1.16e-01} & 0.90\\
3 & \num{2.05e-01} & \num{9.43e-02} & \uline{\num{6.70e-02}} & \num{6.23e-02} & \num{6.14e-02} & \num{6.14e-02} & 0.92\\
4 & \num{1.82e-01} & \num{6.99e-02} & \num{3.83e-02} & \uline{\num{3.28e-02}} & \num{3.19e-02} & \num{3.18e-02} & 0.95\\
5 & \num{1.71e-01} & \num{5.82e-02} & \num{2.40e-02} & \num{1.74e-02} & \uline{\num{1.64e-02}} & \num{1.62e-02} & 0.97\\
6 & \num{1.66e-01} & \num{5.34e-02} & \num{1.74e-02} & \num{9.56e-03} & \num{8.43e-03} & \uline{\num{8.24e-03}} & 0.98\\
\midrule
$\eoc{x}$ & -- & 1.63 & 1.62 & 0.87 & 0.18 & 0.03 \\
$\uline{\eoc{xt}}$ & -- & 1.25 & 1.12 & 1.03 & 1.00 & 0.99 \\
\bottomrule
\end{tabular}
\label{tab:kite.bdf1.l2l2}
\end{table}

\begin{table}
\centering
\caption{Convergence results for the BDF1 scheme in the $ L^\infty(L^2)$ norm for the kite geometry test case in \Cref{sec:num-ex.subsec:kite}}
\begin{tabular}{llllllll}
\toprule
$L_t$\textbackslash $L_x$  & 0 & 1 & 2 & 3 & 4 & 5 & $\text{eoc}_{t}$ \\
\midrule
0 & \num{7.73e-01} & \num{5.70e-01} & \num{5.11e-01} & \num{4.95e-01} & \num{4.89e-01} & \num{4.86e-01} & --\\
1 & \uline{\num{4.96e-01}} & \num{3.50e-01} & \num{2.98e-01} & \num{2.84e-01} & \num{2.81e-01} & \num{2.81e-01} & 0.79\\
2 & \num{3.97e-01} & \uline{\num{2.15e-01}} & \num{1.72e-01} & \num{1.60e-01} & \num{1.58e-01} & \num{1.58e-01} & 0.83\\
3 & \num{3.41e-01} & \num{1.51e-01} & \uline{\num{9.90e-02}} & \num{8.80e-02} & \num{8.61e-02} & \num{8.59e-02} & 0.88\\
4 & \num{3.14e-01} & \num{1.19e-01} & \num{6.01e-02} & \uline{\num{4.75e-02}} & \num{4.55e-02} & \num{4.51e-02} & 0.93\\
5 & \num{3.01e-01} & \num{1.04e-01} & \num{4.04e-02} & \num{2.60e-02} & \uline{\num{2.37e-02}} & \num{2.32e-02} & 0.96\\
6 & \num{2.95e-01} & \num{9.68e-02} & \num{3.14e-02} & \num{1.50e-02} & \num{1.23e-02} & \uline{\num{1.19e-02}} & 0.97\\
\midrule
$\text{eoc}_{x}$ & -- & 1.61 & 1.62 & 1.06 & 0.29 & 0.06 \\
$\uline{\text{eoc}_{xt}}$ & -- & 1.21 & 1.12 & 1.06 & 1.01 & 1.00 \\
\bottomrule
\end{tabular}
\label{tab:kite.bdf1.linfl2}
\end{table}

\begin{table}
\centering
\caption{Convergence results for the BDF1 scheme in the $ L^2(H^1)$ norm for the kite geometry test case in \Cref{sec:num-ex.subsec:kite}}
\begin{tabular}{llllllll}
\toprule
$L_t$~\textbackslash~$L_x$  & 0 & 1 & 2 & 3 & 4 & 5 & $\eoc{t}$ \\
\midrule
0 & \num{2.40e+00} & \num{1.93e+00} & \num{1.88e+00} & \num{1.91e+00} & \num{1.94e+00} & \num{1.98e+00} & --\\
1 & \uline{\num{1.71e+00}} & \num{1.40e+00} & \num{1.33e+00} & \num{1.34e+00} & \num{1.37e+00} & \num{1.42e+00} & 0.48\\
2 & \num{1.45e+00} & \uline{\num{1.04e+00}} & \num{8.99e-01} & \num{8.65e-01} & \num{9.01e-01} & \num{9.54e-01} & 0.57\\
3 & \num{1.32e+00} & \num{8.58e-01} & \uline{\num{6.21e-01}} & \num{5.46e-01} & \num{5.65e-01} & \num{6.20e-01} & 0.62\\
4 & \num{1.26e+00} & \num{7.76e-01} & \num{4.77e-01} & \uline{\num{3.59e-01}} & \num{3.45e-01} & \num{3.81e-01} & 0.70\\
5 & \num{1.23e+00} & \num{7.41e-01} & \num{4.16e-01} & \num{2.53e-01} & \uline{\num{2.18e-01}} & \num{2.27e-01} & 0.75\\
6 & \num{1.22e+00} & \num{7.28e-01} & \num{3.93e-01} & \num{2.10e-01} & \num{1.41e-01} & \uline{\num{1.40e-01}} & 0.70\\
\midrule
$\eoc{x}$ & -- & 0.74 & 0.89 & 0.90 & 0.58 & 0.00 \\
$\uline{\eoc{xt}}$ & -- & 0.72 & 0.75 & 0.79 & 0.72 & 0.64 \\
\bottomrule
\end{tabular}
\label{tab:kite.bdf1.l2h1}
\end{table}

\begin{table}
\centering
\caption{Convergence results for the BDF2 scheme in the $ L^2(L^2)$ norm for the kite geometry test case in \Cref{sec:num-ex.subsec:kite}}
\begin{tabular}{llllllll}
\toprule
$L_t$~\textbackslash~$L_x$  & 0 & 1 & 2 & 3 & 4 & 5 & $\eoc{t}$ \\
\midrule
0 & \num{6.20e-01} & \num{4.95e-01} & \num{4.27e-01} & \num{3.85e-01} & \num{3.62e-01} & \num{3.51e-01} & --\\
1 & \uline{\num{4.16e-01}} & \num{2.56e-01} & \num{1.90e-01} & \num{1.64e-01} & \num{1.52e-01} & \num{1.46e-01} & 1.27\\
2 & \num{2.45e-01} & \uline{\num{1.33e-01}} & \num{7.80e-02} & \num{6.02e-02} & \num{5.45e-02} & \num{5.22e-02} & 1.48\\
3 & \num{2.01e-01} & \num{7.15e-02} & \uline{\num{3.49e-02}} & \num{2.06e-02} & \num{1.68e-02} & \num{1.59e-02} & 1.71\\
4 & \num{1.81e-01} & \num{5.88e-02} & \num{1.89e-02} & \uline{\num{8.23e-03}} & \num{5.07e-03} & \num{4.46e-03} & 1.84\\
5 & \num{1.71e-01} & \num{5.41e-02} & \num{1.55e-02} & \num{4.41e-03} & \uline{\num{1.85e-03}} & \num{1.25e-03} & 1.84\\
6 & \num{1.65e-01} & \num{5.13e-02} & \num{1.41e-02} & \num{3.57e-03} & \num{9.80e-04} & \uline{\num{4.20e-04}} & 1.57\\
\midrule
$\eoc{x}$ & -- & 1.69 & 1.87 & 1.98 & 1.86 & 1.22 \\
$\uline{\eoc{xt}}$ & -- & 1.65 & 1.92 & 2.09 & 2.15 & 2.14 \\
\bottomrule
\end{tabular}
\label{tab:kite.bdf2.l2l2}
\end{table}

\begin{table}
\centering
\caption{Convergence results for the BDF2 scheme in the $ L^\infty(L^2)$ norm for the kite geometry test case in \Cref{sec:num-ex.subsec:kite}}
\begin{tabular}{llllllll}
\toprule
$L_t$\textbackslash $L_x$  & 0 & 1 & 2 & 3 & 4 & 5 & $\text{eoc}_{t}$ \\
\midrule
0 & \num{7.52e-01} & \num{6.04e-01} & \num{5.10e-01} & \num{4.43e-01} & \num{4.01e-01} & \num{3.79e-01} & --\\
1 & \uline{\num{5.79e-01}} & \num{3.54e-01} & \num{2.53e-01} & \num{2.09e-01} & \num{1.88e-01} & \num{1.77e-01} & 1.10\\
2 & \num{3.78e-01} & \uline{\num{1.99e-01}} & \num{1.09e-01} & \num{7.79e-02} & \num{6.80e-02} & \num{6.44e-02} & 1.46\\
3 & \num{3.38e-01} & \num{1.16e-01} & \uline{\num{5.36e-02}} & \num{2.81e-02} & \num{2.09e-02} & \num{1.93e-02} & 1.74\\
4 & \num{3.10e-01} & \num{1.01e-01} & \num{3.29e-02} & \uline{\num{1.29e-02}} & \num{6.65e-03} & \num{6.22e-03} & 1.64\\
5 & \num{2.99e-01} & \num{9.60e-02} & \num{2.83e-02} & \num{7.79e-03} & \uline{\num{2.83e-03}} & \num{1.96e-03} & 1.67\\
6 & \num{2.94e-01} & \num{9.25e-02} & \num{2.60e-02} & \num{6.57e-03} & \num{1.73e-03} & \uline{\num{6.12e-04}} & 1.68\\
\midrule
$\text{eoc}_{x}$ & -- & 1.67 & 1.83 & 1.98 & 1.93 & 1.50 \\
$\uline{\text{eoc}_{xt}}$ & -- & 1.54 & 1.89 & 2.06 & 2.19 & 2.21 \\
\bottomrule
\end{tabular}
\label{tab:kite.bdf2.linfl2}
\end{table}

\begin{table}
\centering
\caption{Convergence results for the BDF2 scheme in the $ L^2(H^1)$ norm for the kite geometry test case in \Cref{sec:num-ex.subsec:kite}}
\begin{tabular}{llllllll}
\toprule
$L_t$~\textbackslash~$L_x$  & 0 & 1 & 2 & 3 & 4 & 5 & $\eoc{t}$ \\
\midrule
0 & \num{2.62e+00} & \num{2.20e+00} & \num{2.03e+00} & \num{1.96e+00} & \num{1.97e+00} & \num{2.02e+00} & --\\
1 & \num{2.02e+00} & \num{1.45e+00} & \num{1.23e+00} & \num{1.17e+00} & \num{1.17e+00} & \num{1.22e+00} & 0.73\\
2 & \dashuline{\num{1.46e+00}} & \num{1.02e+00} & \num{7.66e-01} & \num{6.42e-01} & \num{6.06e-01} & \num{6.17e-01} & 0.98\\
3 & \num{1.33e+00} & \uwave{\num{8.04e-01}} & \dashuline{\num{5.46e-01}} & \num{3.69e-01} & \num{2.98e-01} & \num{2.92e-01} & 1.08\\
4 & \num{1.27e+00} & \num{7.59e-01} & \uuline{\num{4.28e-01}} & \uwave{\num{2.56e-01}} & \num{1.64e-01} & \num{1.43e-01} & 1.03\\
5 & \num{1.24e+00} & \num{7.38e-01} & \num{4.03e-01} & \uuline{\num{2.06e-01}} & \uline{\num{1.16e-01}} & \num{7.68e-02} & 0.90\\
6 & \num{1.22e+00} & \num{7.26e-01} & \num{3.90e-01} & \num{1.95e-01} & \num{9.74e-02} & \uuline{\num{5.37e-02}} & 0.52\\
\midrule
$\eoc{x}$ & -- & 0.75 & 0.90 & 1.00 & 1.00 & 0.86 \\
$\eoc{xxt}$ & -- & -- & 1.42 & 1.65 & 1.89 & 1.94 \\
\bottomrule
\end{tabular}
\label{tab:kite.bdf2.l2h1}
\end{table}

\subsection{Example 3: Colliding circles}
\label{sec:num-ex.subsec:colliding}

This more advanced example consists of two circles that collide and separate again~\cite{LO19}. Consequently, this includes a topology change in the geometry and a discontinuous transport field. Here we do not have an analytical solution, However, we can track the conservation of our scalar quantity.

The geometry is described by the level set function
\begin{equation*}
  \phi(\xb, t) = \min\{\Vert\xb - s_1(t)\Vert_2, \Vert\xb - s_1(t)\Vert_2\} - R,\text{ with }s_1(t) = (0, t-3/4)^T,\; s_1(t) = (0, 3/4 - t)^T. 
\end{equation*}
The radius is chosen as $R=0.5$ and $\tend=1.5$, such that $\phi(\xb, 0)=\phi(\xb, \tend)$. The transport field is given by
\begin{equation*}
  \wb = 
  \begin{cases}
    (0, -1)^T & \text{if ($\xb_2>0$ and $t\leq \tend/2$) or ($\xb_2\leq0$ and $t> \tend/2$)}\\
    (0, 1)^T  & \text{if ($\xb_2\leq0$ and $t\leq \tend/2$) or ($\xb_2>0$ and $t< \tend/2$)}.
  \end{cases}
\end{equation*}
The background domain considered is $\widetilde{\O} = (-0.6, 0.6)\times(-1.35, 1.35)$ and diffusion coefficient is chosen to be $0.1$. Finally, the initial condition is given as $u_0 = \sign(\xb_2)$.

\begin{figure}
  \centering
  -1 \includegraphics[height=6.2pt, width=0.8\textwidth]{img/colliding_h0.07dt0.01875colourbar.png} 1\phantom{-}\\[4pt]
  \includegraphics[width=0.089\textwidth]{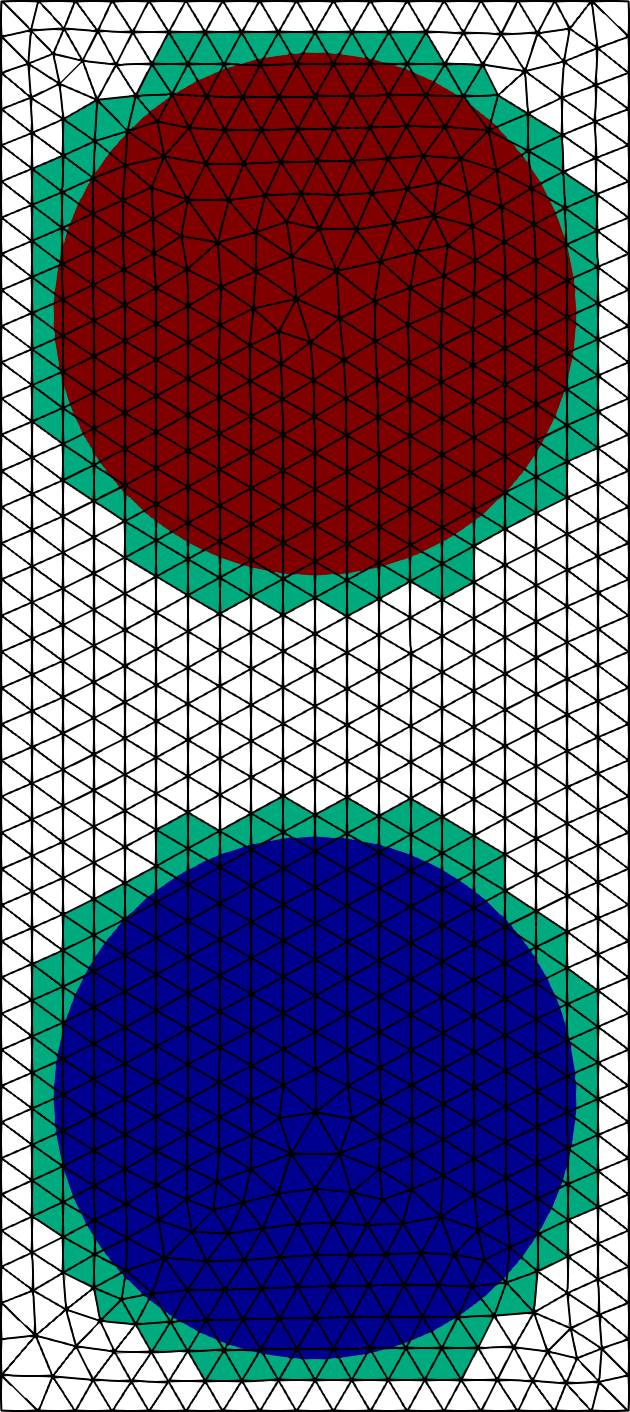}\hspace{.5pt}%
  \includegraphics[width=0.089\textwidth]{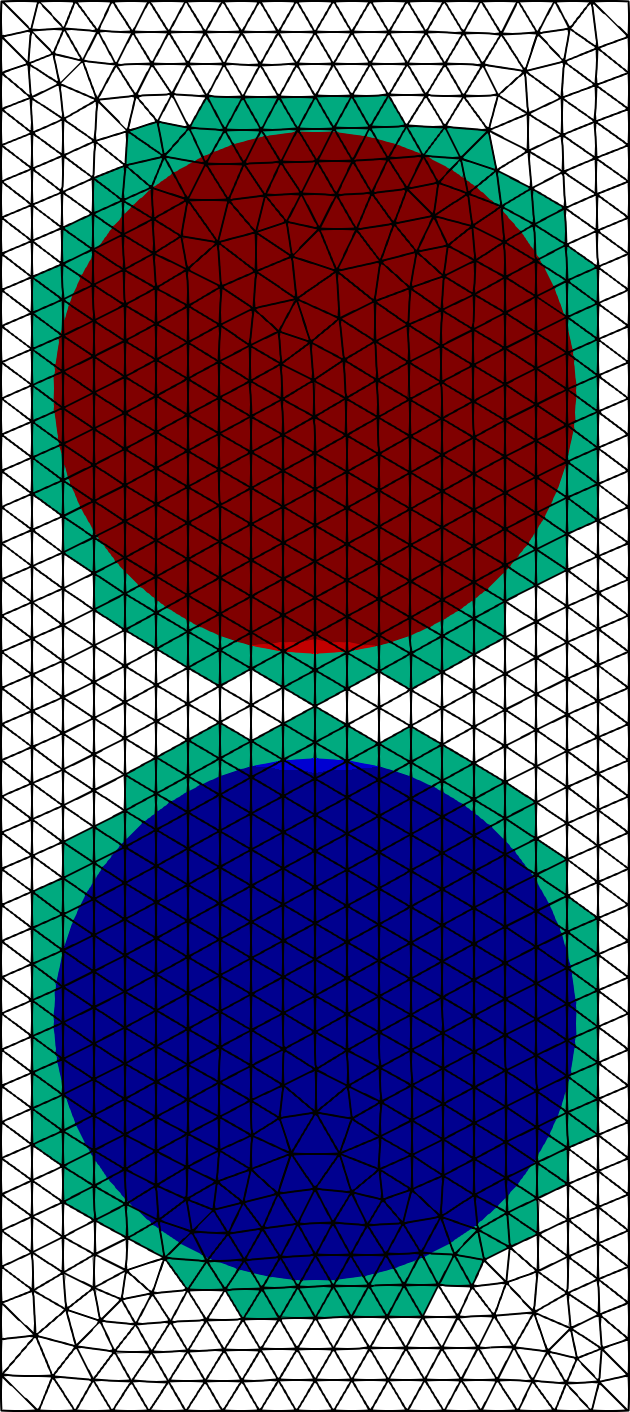}\hspace{.5pt}%
  \includegraphics[width=0.089\textwidth]{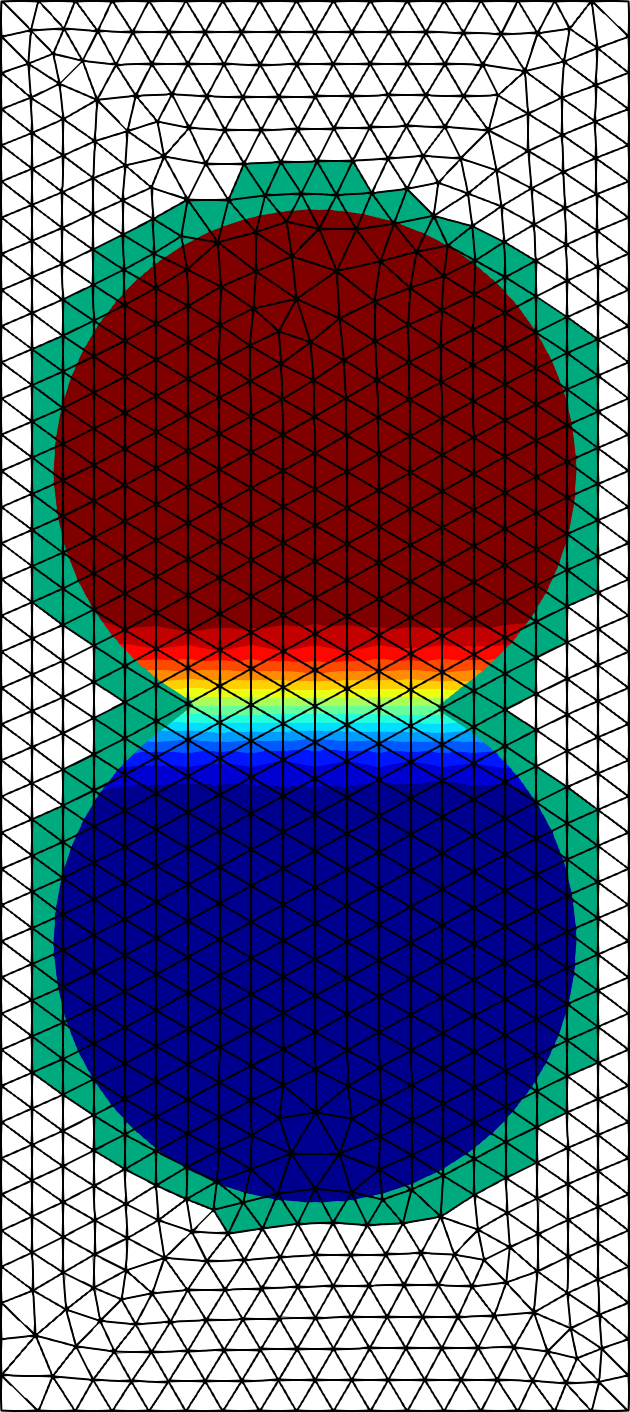}\hspace{.5pt}%
  \includegraphics[width=0.089\textwidth]{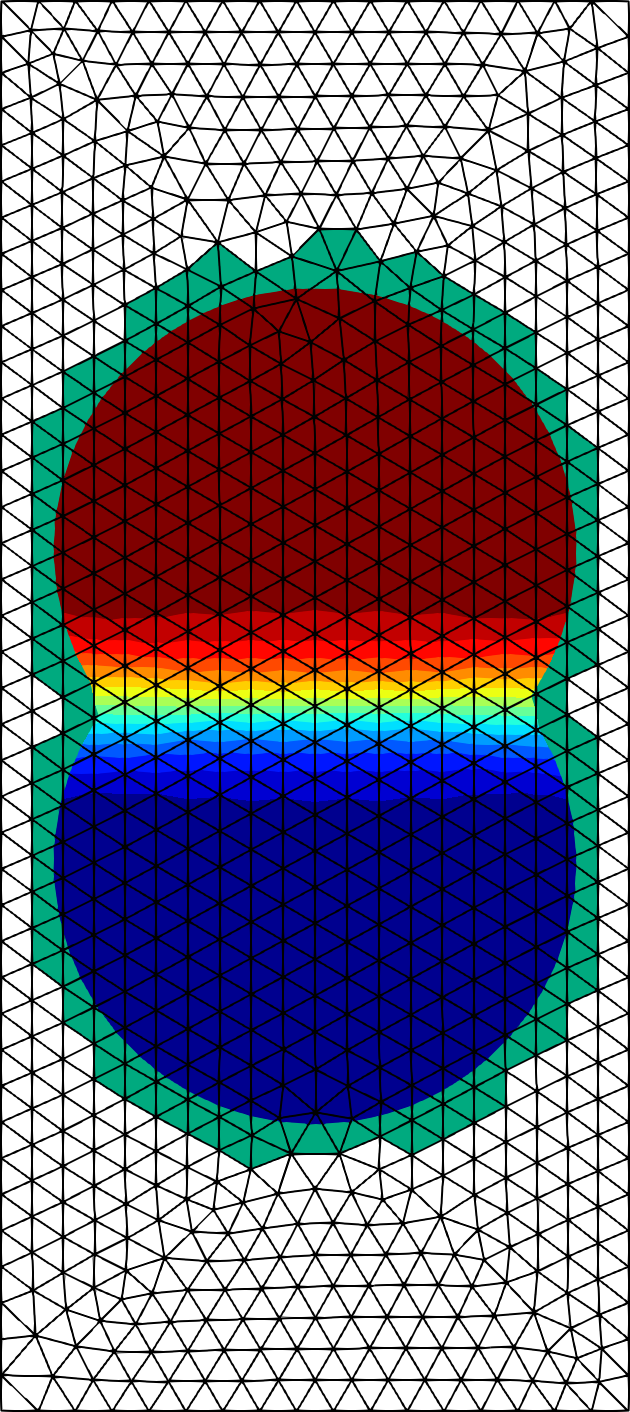}\hspace{.5pt}%
  \includegraphics[width=0.089\textwidth]{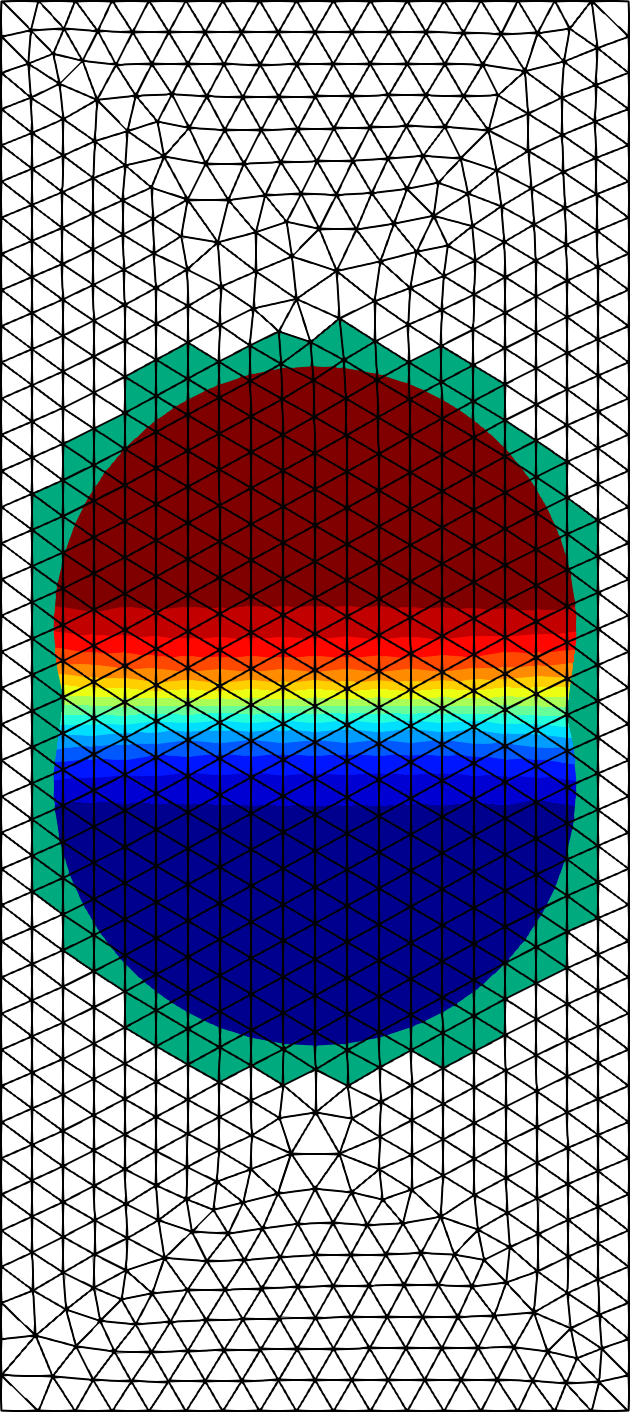}\hspace{.5pt}%
  \includegraphics[width=0.089\textwidth]{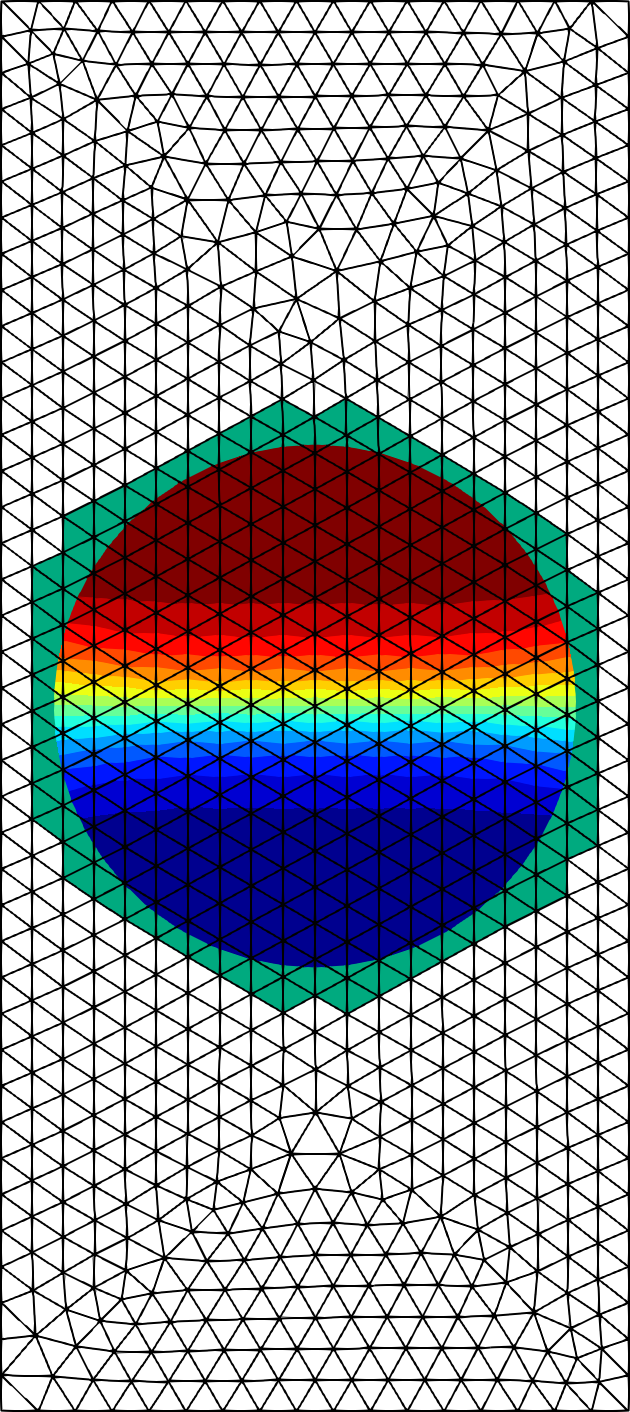}\hspace{.5pt}%
  \includegraphics[width=0.089\textwidth]{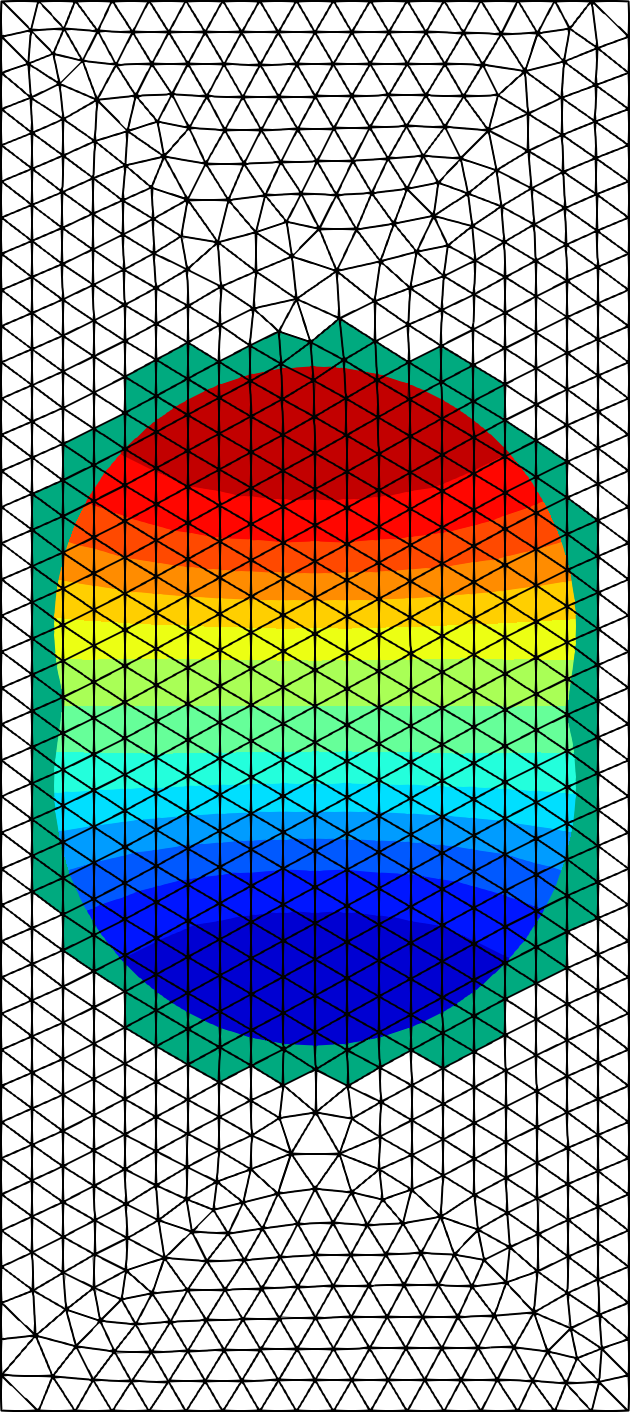}\hspace{.5pt}%
  \includegraphics[width=0.089\textwidth]{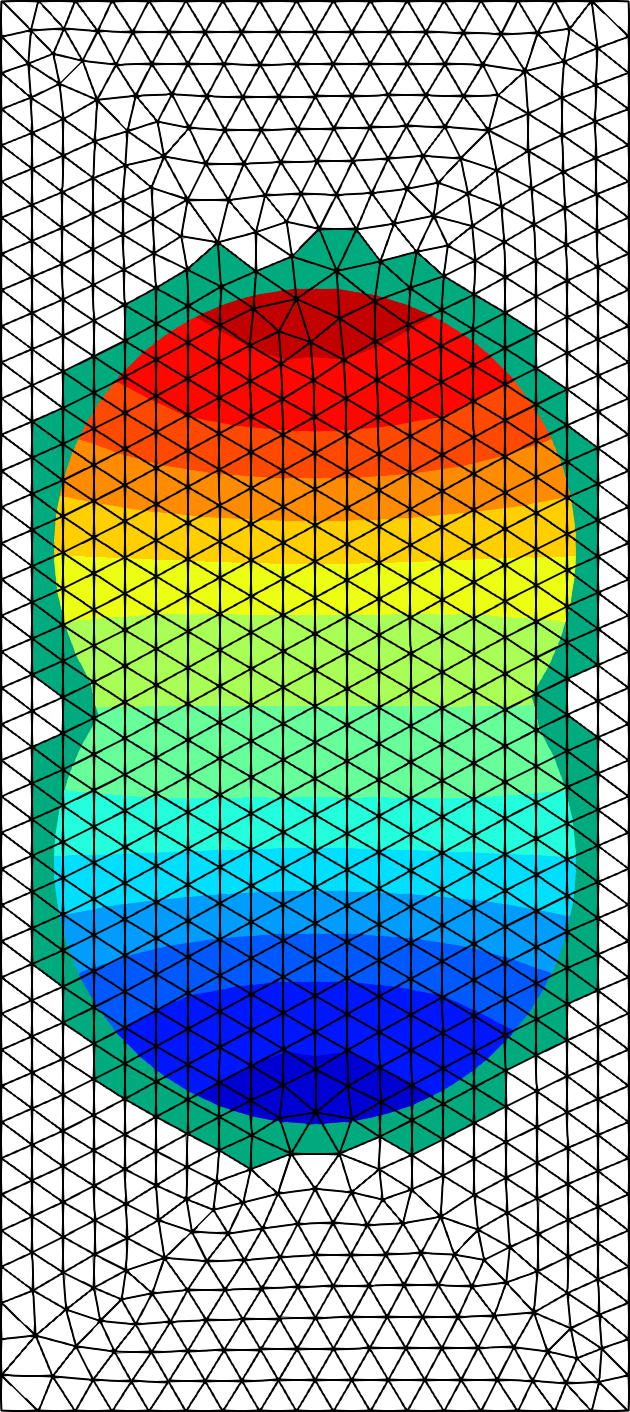}\hspace{.5pt}%
  \includegraphics[width=0.089\textwidth]{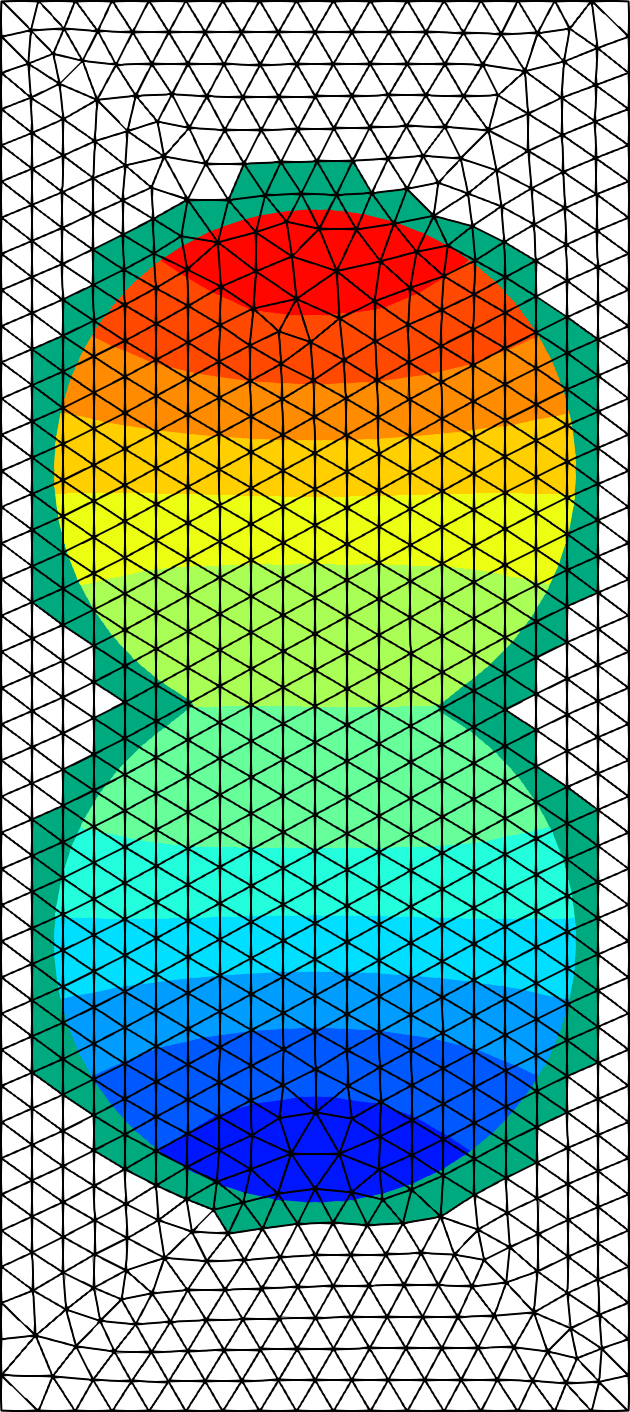}\hspace{.5pt}%
  \includegraphics[width=0.089\textwidth]{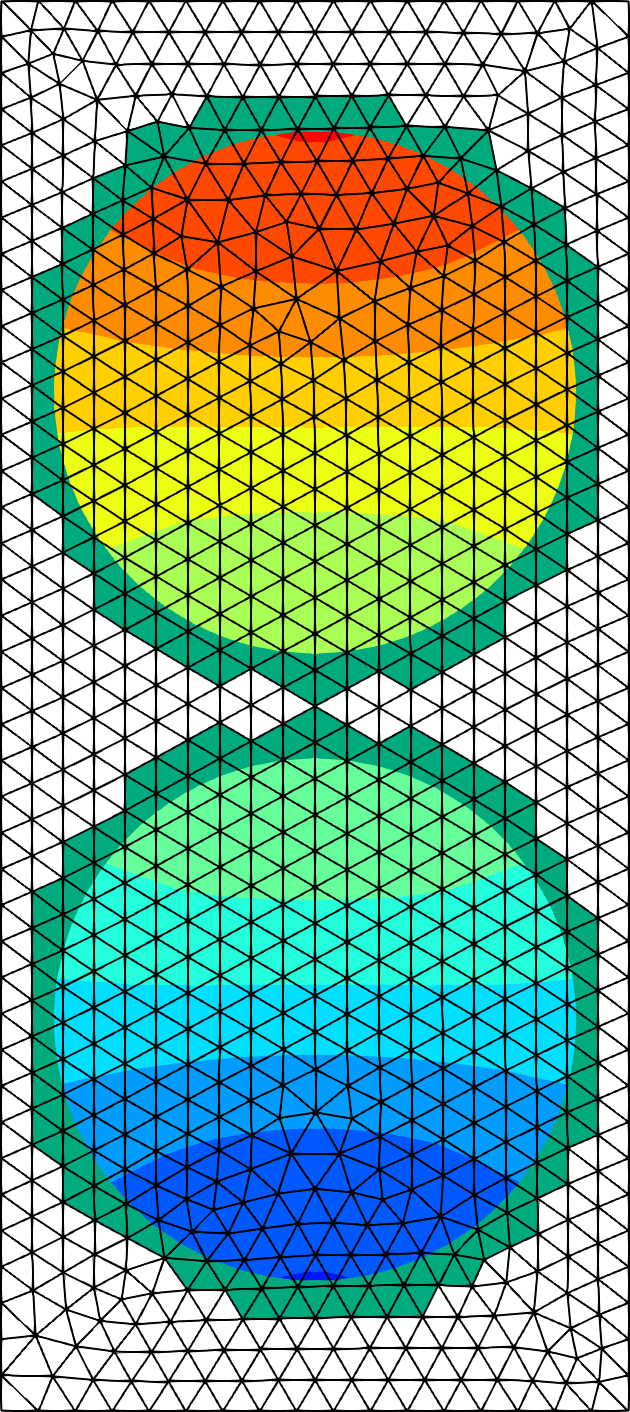}\hspace{.5pt}%
  \includegraphics[width=0.089\textwidth]{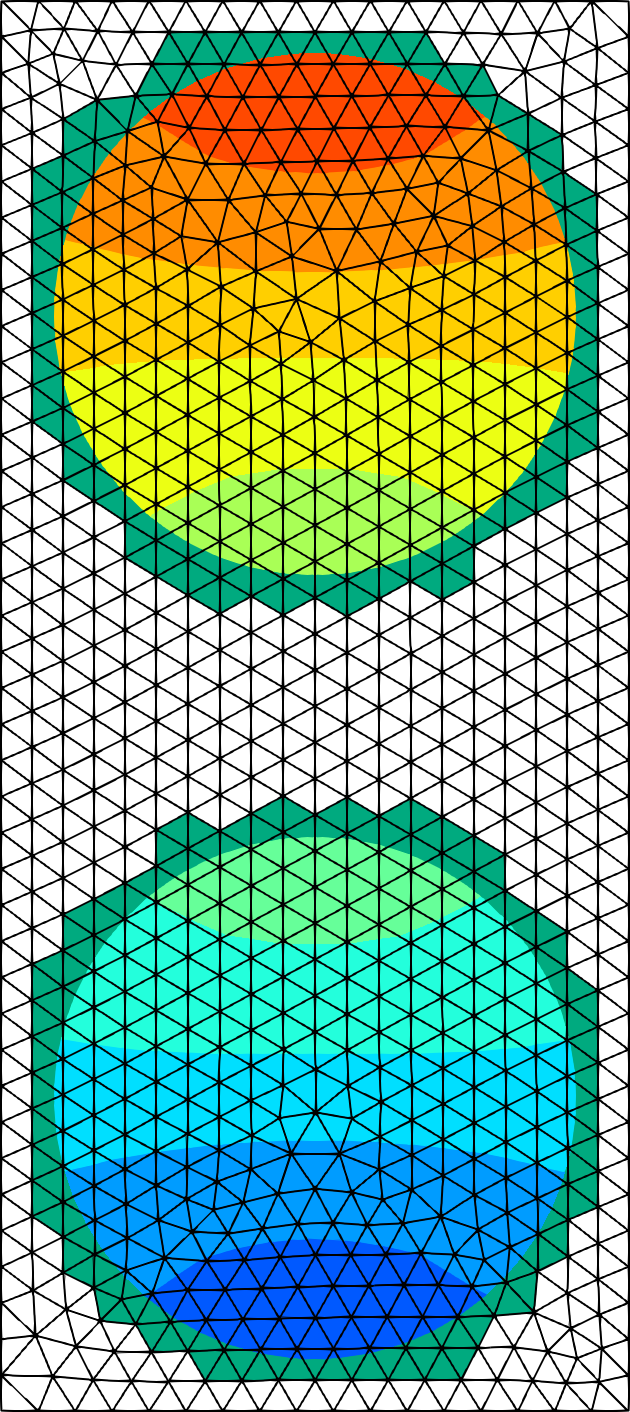}
  \caption{Discrete solution for the colliding circles test case in \Cref{sec:num-ex.subsec:colliding} at intervals of $0.1\tend$. Extension elements in $\ThS$ are marked in green.}
  \label{fig.colliding-cir}
\end{figure}

We take $h=0.07$, $\dt=T/80$ and use the BDF2 scheme. The results at intervals of $0.1T$ can be seen in \Cref{fig.colliding-cir}. Visually, these results match those presented in~\cite{LO19}; however, here we conserve the total of our scalar quantity up to machine precision in every time step. As in~\cite{LO19}, mass is exchanged between the two domains, as soon as the ghost-penalty extension domains overlap, rather than when the domains overlap.

\subsection{Example 4: Colliding spheres}
\label{sec:num-ex.subsec:colliding3d}

Finally, we extend the previous example to three spatial dimensions. The geometry is now described by the level set function
\begin{equation*}
  \phi(\xb, t) = \min\{\Vert\xb - s_1(t)\Vert_2, \Vert\xb - s_1(t)\Vert_2\} - R,\text{ with }s_1(t) = (0, 0, t-3/4)^T,\; s_1(t) = (0, 0, 3/4 - t)^T. 
\end{equation*}
The radius is chosen as $R=0.5$, the end time as $\tend=1.5$, and the transport field as
\begin{equation*}
  \wb = 
  \begin{cases}
    (0, 0, -1)^T & \text{if ($\xb_3>0$ and $t\leq \tend/2$) or ($\xb_3\leq0$ and $t> \tend/2$)}\\
    (0, 0, 1)^T  & \text{if ($\xb_3\leq0$ and $t\leq \tend/2$) or ($\xb_3>0$ and $t< \tend/2$)}.
  \end{cases}
\end{equation*}
The background domain is $\widetilde{\O} = (-0.6, 0.6)\times(-0.6, 0.6)\times(-1.35, 1.35)$, $\nu=0.1$ and  $u_0 = \sign(\xb_3)$.

\begin{figure}
  \centering
  -1 \includegraphics[height=6.2pt, width=0.8\textwidth]{img/colliding_h0.07dt0.01875colourbar.png} 1\phantom{-}\\[4pt]
  \includegraphics[width=0.15\textwidth]{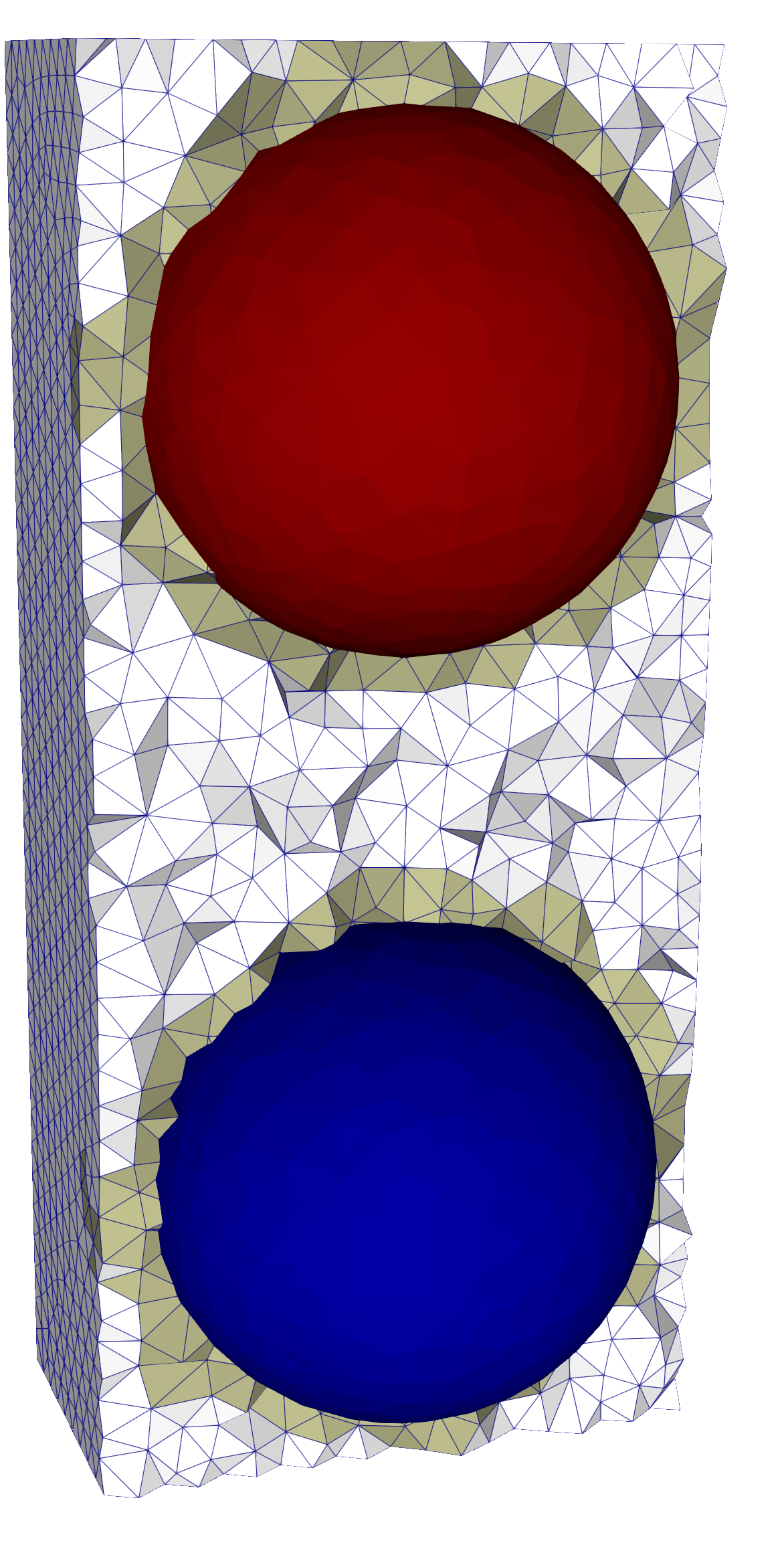}\hspace{.5pt}%
  \includegraphics[width=0.15\textwidth]{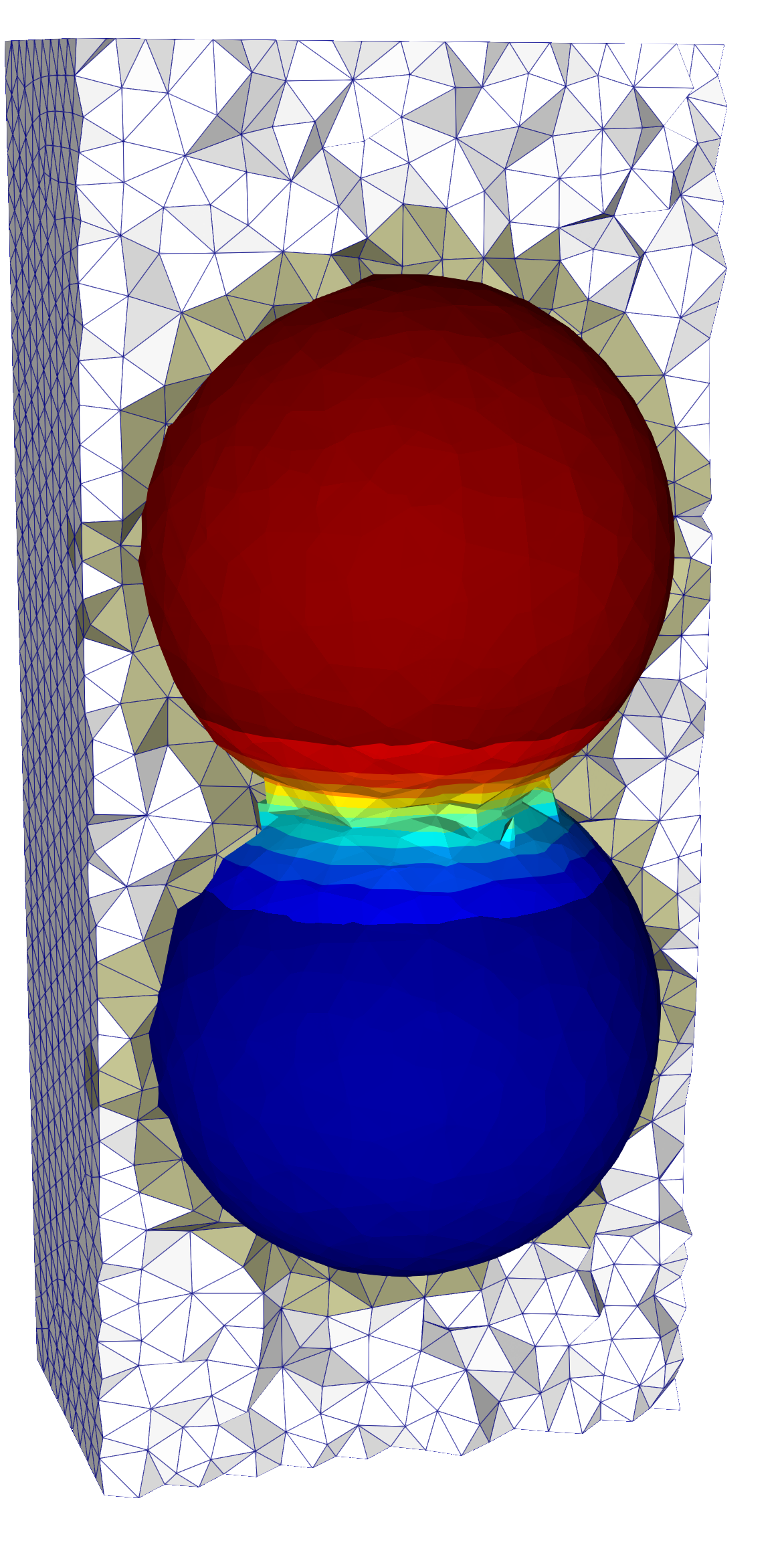}\hspace{.5pt}%
  \includegraphics[width=0.15\textwidth]{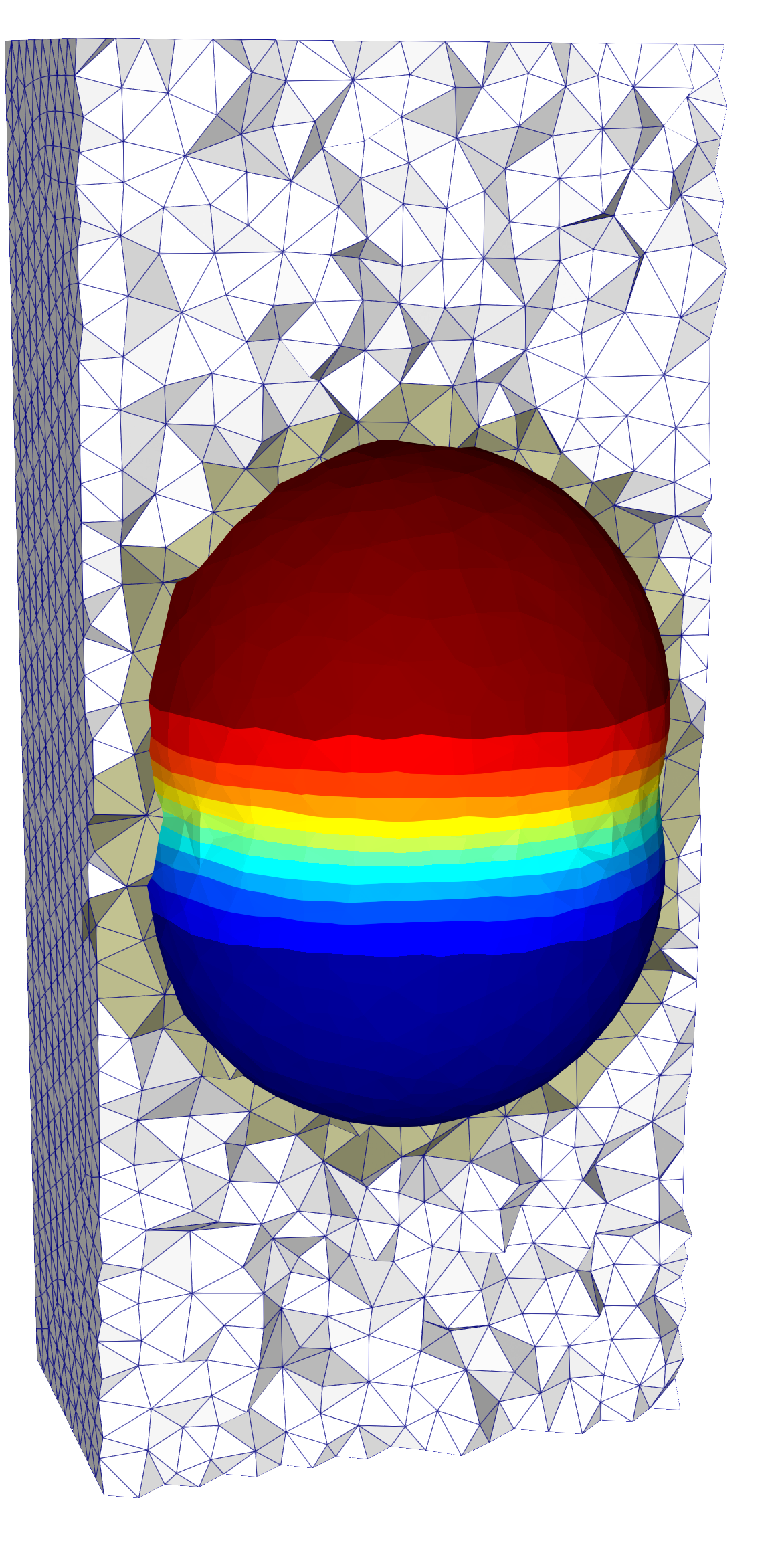}\hspace{.5pt}%
  \includegraphics[width=0.15\textwidth]{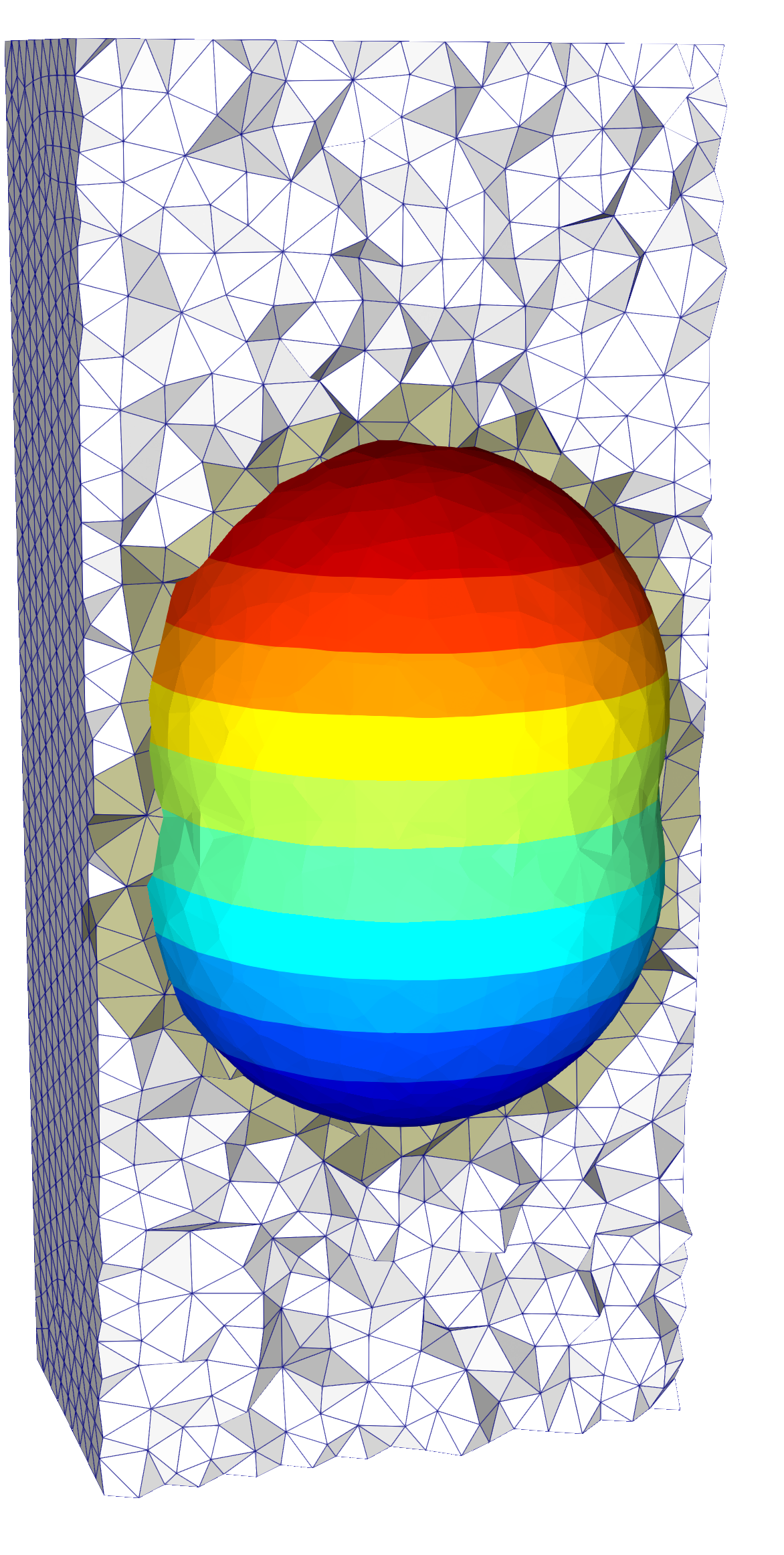}\hspace{.5pt}%
  \includegraphics[width=0.15\textwidth]{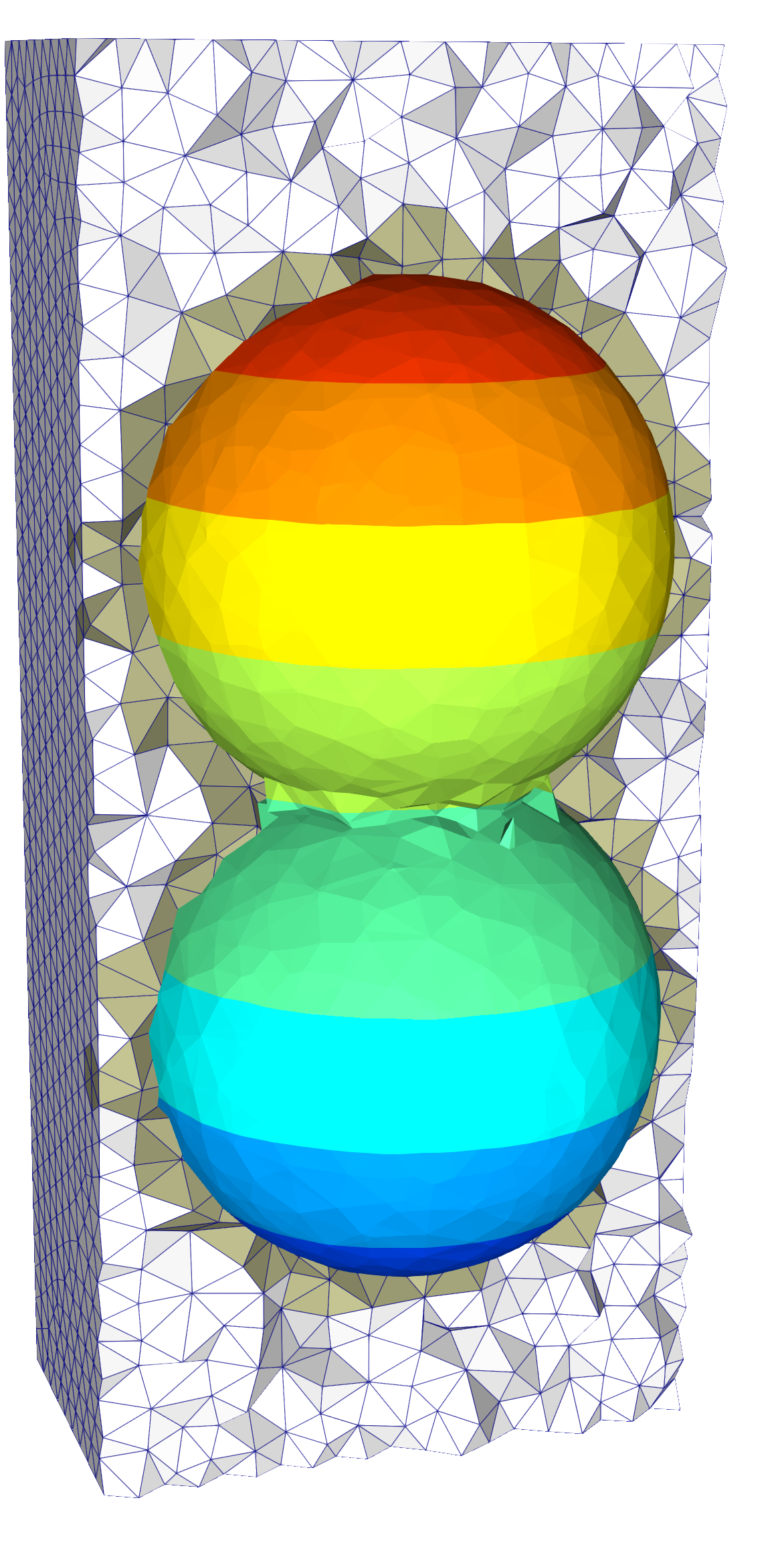}\hspace{.5pt}%
  \includegraphics[width=0.15\textwidth]{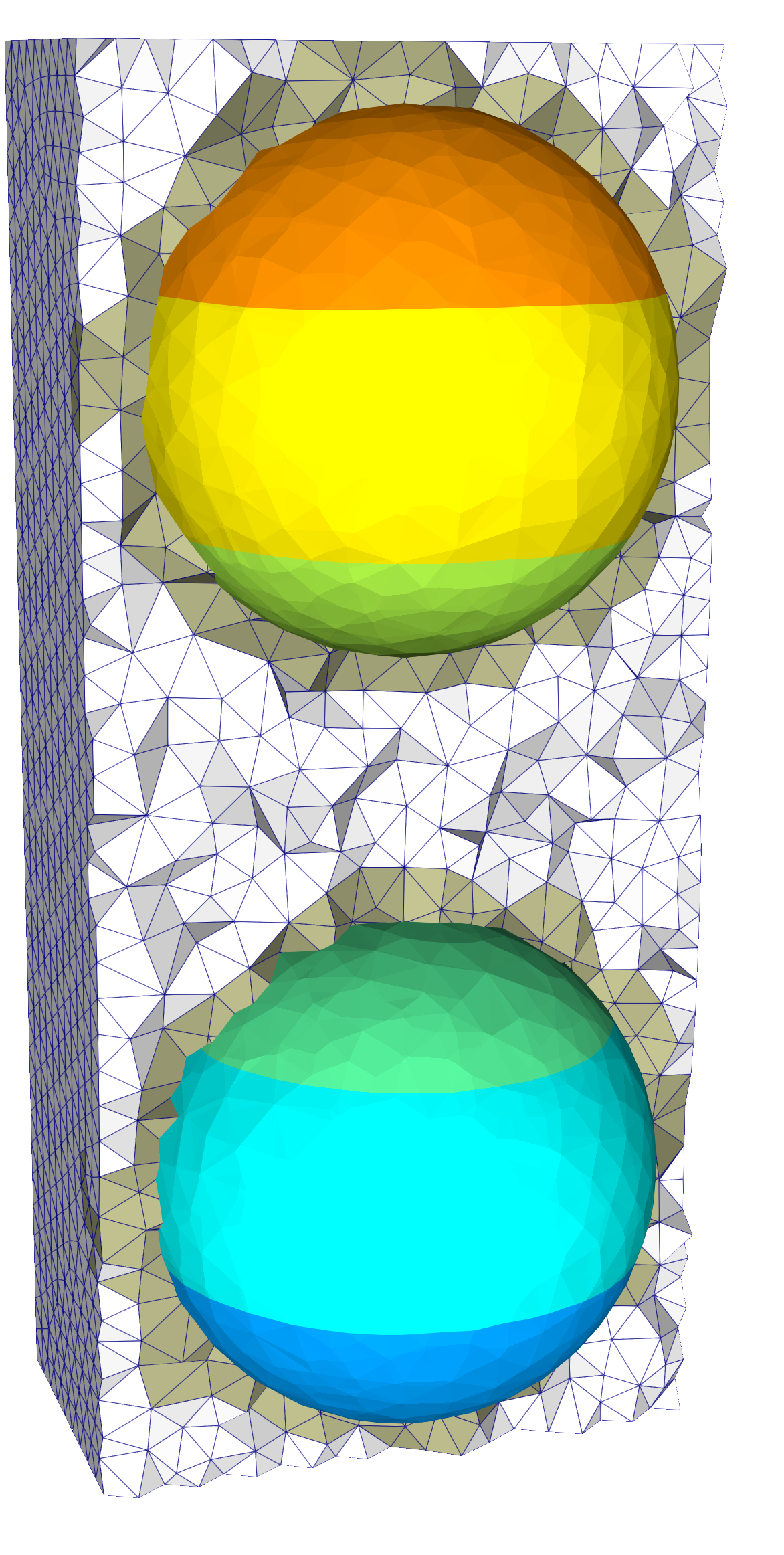}
  \caption{Discrete solution for the colliding sphere test case in \Cref{sec:num-ex.subsec:colliding3d} at intervals of $0.2T$. Extension elements in $\ThS$ are marked in green.}
  \label{fig.colliding-sphere}
\end{figure}

We again take $h=0.07$, $\dt=T/80$ and use the BDF2 version of our time-stepping scheme. The results at intervals of $0.2\tend$ can be seen in \Cref{fig.colliding-sphere}. We again observe the same behavior as in the previous two-dimensional behavior, and the total of the scalar quantity is preserved up to machine precision in every time step.

\section*{Data Availability Statement}
The code used to realise the numerical examples is freely available on github \url{https://github.com/hvonwah/conserv-eulerian-moving-domian-repro} and archived on zenodo \url{https://doi.org/10.5281/zenodo.10951768}.

\section*{Acknowledgements}
This material is based upon work supported by the National Science Foundation
under Grant No. DMS-1929284 while the authors were in residence at the
Institute for Computational and Experimental Research in Mathematics in
Providence, RI, during the Numerical PDEs: Analysis, Algorithms, and Data
Challenges program. The author M.O. was partially supported by the National 
Science Foundation grants DMS-2309197 and DMS-2408978.

\printbibliography

\end{document}